\documentclass[a4paper,12pt]{amsart}
\usepackage{amssymb,amsmath,a4wide,graphicx,stmaryrd,fullpage,setspace,microtype,textcomp,tikz,tikz-cd,accents,mathtools,abraces}
\usepackage{hyperref}
\usepackage[shortlabels]{enumitem}
\definecolor{dark-red}{rgb}{0.5,0.15,0.15}
\hypersetup{
	colorlinks   = true, %Colours links instead of ugly boxes
	urlcolor     = dark-red, %Colour for external hyperlinks
	linkcolor    = black, %Colour of internal links
	citecolor   = dark-red %Colour of citations
}
\usepackage[T1]{fontenc}
\usepackage{lmodern}
\usepackage[numbers,sort&compress]{natbib}

\title{Globular subdivisions are dihomotopy equivalences}

\author[P. Gaucher]{Philippe Gaucher}

\address{Universit\'e Paris Cit\'e, CNRS, IRIF, F-75013, Paris, France}

\urladdr{http://www.irif.fr/{\~{}}gaucher} 

\makeatletter
\@namedef{subjclassname@2020}{%
	\textup{2020} Mathematics Subject Classification}
\makeatother


\swapnumbers

\makeatletter
\let\P\@undefined
\let\leq\@undefined
\let\geq\@undefined
\let\top\@undefined
\makeatother
\newcommand{\C}{\mathcal{C}}
\newcommand{\D}{\mathcal{D}}
\newcommand{\W}{\mathcal{W}}
\newcommand{\F}{\mathcal{F}}
\newcommand{\K}{\mathcal{K}}
\newcommand{\T}{\mathcal{T}}
\newcommand{\de}{\partial}
\newcommand{\p}{\times}
\renewcommand{\vec}{\overrightarrow}
\newcommand{\P}{\mathbb{P}}
\DeclareMathOperator{\II}{I^+}
\newcommand{\cont}{{\vec{\mathrm{Sp}}}}

\newcommand{\discont}{{\vec{\Omega}}}
\newcommand{\dtop}{{\brm{Flow}}}
\newcommand{\dcat}{{\mathrm{cat}}}
\newcommand{\moore}{{\mathbb{M}}}
\newcommand{\lmoore}{\mathbb{M}_!}
\newcommand{\glob}{{\rm{Glob}}}
\DeclareMathOperator{\cell}{{\brm{cell}}}
\DeclareMathOperator{\cof}{{\brm{cof}}}
\DeclareMathOperator{\inj}{{\brm{inj}}}
\newcommand{\rest}{{\!\restriction}}
\newcommand{\ddownarrow}{{\downarrow}}
\newcommand{\sbd}{\mathrm{sbd}}
\newcommand{\top}{{\mathbf{Top}}}
\newcommand{\mtop}{\mathbf{mTop}}
\newcommand{\iso}{\cong}
\newcommand{\vI}{\vec{I}}
\newcommand{\leq}{\leqslant}
\newcommand{\geq}{\geqslant}
\newcommand{\tr}[1]{{\langle{#1}\rangle}}
\newcommand{\brm}[1]{\rm{\mathbf{#1}}}
\newcommand{\ttop}{{\brm{TOP}}}
\newcommand{\globM}{{\rm{Glob}}^{top}}
\DeclareMathOperator{\id}{Id}
\newcommand{\liminj}{\varinjlim}
\newcommand{\limproj}{\varprojlim}
\newcommand{\cat}{{\mathbf{Cat}}}
\newcommand{\ptop}[1]{{\brm{{#1}dTop}}}
\newcommand{\cocartesian}{\arrow[lu, phantom, "\ulcorner"{font=\Large}, pos=0]}

\newcommand{\ttt}{two-out-of-three property}
\newtheorem*{thmN}{Theorem}
\newtheorem{thm}{Theorem}[section]
\newtheorem{prop}[thm]{Proposition}

\newtheorem{cor}[thm]{Corollary}
\newcommand{\bth}{\begin{thm}}
\renewcommand{\eth}{\end{thm}}
\newcommand{\bpf}{\begin{proof}}
\newcommand{\epf}{\end{proof}}
\theoremstyle{definition}
\newtheorem{defn}[thm]{Definition}
\newcommand{\bd}{\begin{defn}}
\newcommand{\ed}{\end{defn}}
\newtheorem{nota}[thm]{Notation}
\makeatletter
\newcommand*{\@opargbegintheorem}[3]{\trivlist
	\item[\hskip \labelsep{\bfseries #1\ #2}] \textbf{(#3)}\ \itshape}
\makeatother

\allowdisplaybreaks

\begin{document}

\begin{abstract} 
	We prove that any globular subdivision of multipointed $d$-spaces gives rise to a dihomotopy equivalence between the associated flows. As a straightforward application, the flows associated to two multipointed $d$-spaces related by a finite zigzag of globular subdivisions have isomorphic branching and merging homology theories and isomorphic underlying homotopy types. 
\end{abstract}

\maketitle
\tableofcontents
\hypersetup{linkcolor = dark-red}

\section{Introduction}

\subsection*{Presentation}

Directed Algebraic Topology (DAT) studies mathematical objects arising from the geometric study of concurrent systems up to homotopy \cite{DAT_book}. There are two main classes of geometric models: the \textit{continuous} models like Grandis' directed spaces (Definition~\ref{def:directed_space}) and Krishnan's streams \cite{MR2545830}, and the \textit{multipointed} models like multipointed $d$-spaces (Definition~\ref{def:multipointed-d-space}) and flows (Definition~\ref{def:flow}) which have a distinguished set of states. The main multipointed \textit{combinatorial} model of concurrency is the category of precubical sets, a $n$-cube representing the concurrent execution of $n$ actions \cite{DAT_book}. A precubical set can be realized in any of the geometric models of concurrency above.

There is no known convenient model category structure on continuous DAT models, as all attempts to date cannot prevent the directed segment from being contracted by weak equivalences, which implies that the weak equivalences destroy the causal structure. 

On the other hand, combinatorial model category structures have been introduced for multipointed $d$-spaces and flows such that the weak equivalences do not contract the directed segment. They are related by a zigzag of Quillen equivalences by \cite[Theorem~10.9]{Moore1} and \cite[Theorem~14]{Moore3}. There is also a functor $\dcat:\ptop{\mathcal{M}} \to \dtop$ from multipointed $d$-spaces to flows which is neither a left adjoint nor a right adjoint and such that the total left derived functor in the sense of \cite{HomotopicalCategory} induces an equivalence of categories between the homotopy categories of the model structures \cite[Theorem~15]{Moore3}. Unfortunately, there are two issues: 1) the weak equivalences are extremely rigid, the weak equivalences inducing bijections between the distinguished sets of states, which implies that a directed segment cannot be identified with the concatenation of two directed segments by a map preserving the extremities; 2) the weak equivalences do not preserve the causal structure \textit{between} the distinguished set of states. The latter problem is explained in detail in \cite[Section~10 and Section~11]{GlobularNaturalSystem}. 

To overcome the first issue, the paper \cite{hocont} introduces a notion of flow up to dihomotopy. Two flows are dihomotopy equivalent if they can be related by a finite zigzag of weak equivalences and of retracts of transfinite compositions of pushouts of generating subdivisions in the sense of Definition~\ref{def:generating-sbd}. A Whitehead theorem is even proved \cite[Theorem~4.6 and 4.7]{hocont}, namely that any dihomotopy equivalence between two q-cofibrant homotopy continuous flows is invertible up to homotopy. The notion of homotopy continuous object plays the role of fibrant object in this theory. Unfortunately, it can also be proved that the categorical localization of flows up to dihomotopy is not the homotopy category of a model category structure on flows \cite[Theorem~5.7]{nonexistence}. Not much more is known about this categorical localization. However, it is conjecturally the homotopy category of a cofibration category, which would be an interesting setting for doing homotopy theory. 

To overcome the second issue, the paper \cite{GlobularNaturalSystem} proves the globular analogue of Dubut's results. More precisely, the latter paper proves that the globular subdivisions (Definition~\ref{def:globular-sbd}) preserve the causal structure between the distinguished sets of states in the following sense. The natural systems in Dubut's sense of the associated directed spaces of two cellular multipointed $d$-spaces related by a finite zigzag of globular subdivisions are bisimilar up to homotopy, whether the discrete or the continuous versions of natural system is used \cite[Theorem~8.16 and Theorem~9.4]{GlobularNaturalSystem}. In particular, this opens the way to computational techniques by representing a directed path as a finite sequence of globular cells. Since the full subcategory of cellular multipointed $d$-spaces (i.e. the cellular objects of their q-model structure) contains all examples coming from computer science, their categorical localization by the globular subdivisions deserves to be studied more carefully. 

This categorical localization, which can be proved to be locally small, has no known mathematical structure. As a first step to understand it, this paper links the approaches of \cite{hocont} and \cite{GlobularNaturalSystem} by the theorem stated now:

\begin{thmN} (Corollary~\ref{final})
	Let $X$ and $Y$ be two cellular multipointed $d$-spaces related by a finite zigzag sequence of globular subdivisions. Then the associated flows $\dcat(X)$ and $\dcat(Y)$ are related by a finite zigzag of transfinite compositions of pushouts of generating subdivisions and of weak equivalences of flows. Using the language of \cite{hocont}, the two flows  $\dcat(X)$ and $\dcat(Y)$ are dihomotopy equivalent. Moreover, the generating subdivisions are of the form a q-cofibrant replacement of the inclusion of the poset $\{0<1\}$ into a globular poset as depicted in Figure~\ref{fig:glposet}.
\end{thmN}

The induced functor from the category of cellular multipointed $d$-spaces up to globular subdivisions to the category of flows up to dihomotopy equivalences cannot be an equivalence of categories. Indeed, \cite[Section~10]{GlobularNaturalSystem} provides in \cite[Proposition~10.1]{GlobularNaturalSystem} an example of a trivial q-fibration between two cellular multipointed $d$-spaces $f:A\to B$ such that the natural systems in Dubut's sense of the associated directed spaces $\cont(A)$ and $\cont(B)$ (see Notation~\ref{notation:cont}) cannot be bisimilar up to homotopy. By \cite[Theorem~8.16 and Theorem~9.4]{GlobularNaturalSystem}, this implies that the cellular multipointed $d$-spaces $A$ and $B$ cannot be isomorphic in the categorical localization of the cellular multipointed $d$-spaces by the globular subdivisions.

It would be interesting to prove the same results for the cubical subdivisions, since cubical subdivisions of realizations of precubical sets as directed spaces preserve the bisimilarity type \cite{dubut_PhD}. This problem is left for a future work.

As a first application of Corollary~\ref{final}, if $f:X\to Y$ is a globular subdivision, then the associated map of flows $\dcat(f):\dcat(X)\to \dcat(Y)$ induces an isomorphism on the branching and merging homology theories thanks to \cite[Corollary~11.3]{3eme}. These homology theories detect the non-deterministic branching and merging areas of execution paths. There are methods to define the branching and merging homology theories of a multipointed d-space without using the functor $\dcat:\ptop{\mathcal{M}}\to \dtop$ and the complicated homotopical constructions of \cite{3eme}. These methods will make the preservation of the branching and merging homology theories by the globular subdivisions much easier to prove. This will be the subject of a subsequent paper.

Another application of Corollary~\ref{final} is that, for any globular subdivision $f:X\to Y$, the associated map of flows $\dcat(f):\dcat(X)\to \dcat(Y)$ induces an isomorphism between the underlying homotopy types thanks to \cite[Theorem~9.1]{4eme}. Intuitively, the underlying homotopy type of a flow (see Definition~\ref{def:underlying-type-flow}) is the homotopy type of space obtained after removing the temporal information. We conclude the paper by an alternative proof which uses the fact that a globular subdivision is a homeomorphism and without using the complicated homotopical constructions of \cite{4eme}.

\subsection*{Outline of the paper}

Section~\ref{section:cat} recalls some basic facts in category theory and in model category theory. 

Section~\ref{section:space-flow} recalls some basic facts about multipointed $d$-spaces and flows.

Section~\ref{section:dcat} recalls some basic facts about the functor $\dcat:\ptop{\mathcal{M}} \to \dtop$ from multipointed $d$-spaces to flows. The important theorems are Theorem~\ref{thm:rappel-all} and Theorem~\ref{thm:colim}. 

Section~\ref{section:enriched-id} recalls the link between the formalism of multipointed $d$-spaces and flows and the formalism of directed spaces and enriched small categories. This enables us to obtain simpler statements for Theorem~\ref{thm:plus-minus-map-flow} and to use some results from \cite{GlobularNaturalSystem} in this paper.  

Section~\ref{section:cellular} introduces the cellular multipointed $d$-spaces and expounds some preparatory lemmas to study globular subdivisions in Section~\ref{section:gl-sbd}. We introduce in Proposition~\ref{prop:subd-gl-globe} a cellular multipointed $d$-space $\globM(\mathbf{D}^n)_F$ which is the $(n+1)$-dimensional topological globe $\globM(\mathbf{D}^n)$ with a finite set $F$ of additional distinguished states belonging to the interior of the underlying topological space. The main result is Proposition~\ref{prop:raf-cell} which states that the map of multipointed $d$-spaces $\globM(\mathbf{S}^{n-1}) \subset \globM(\mathbf{D}^n)_F$ is cellular. We also introduce the related notion of globular poset: one of them is depicted in Figure~\ref{fig:glposet}. 

Section~\ref{section:gl-sbd} studies the notion of globular subdivision (Definition~\ref{def:globular-sbd}) which was introduced in \cite[Definition~4.10]{diCW} for a very specific kind of cellular multipointed $d$-spaces and in \cite[Definition~9.1]{GlobularNaturalSystem} for general cellular multipointed $d$-spaces. The important notion of connection map of a globular subdivision is introduced in Definition~\ref{def:connectionmaps}. Theorem~\ref{thm:plus-minus-map} is the first key fact of the paper: it is a consequence of the definition of the connection maps and of the specific geometric properties of the topological globes, which are already used in \cite{GlobularNaturalSystem}. Finally, Theorem~\ref{thm:sbd-cell} proves that the source and target of any globular subdivision have compatible cellular decompositions. It is the second key fact of the paper. 

Section~\ref{section:glT} introduces the notion of generating subdivision and also a specific subset $\T^{gl}$ of them (see Notation~\ref{notation:Igl}) which is specially designed for the calculations of Proposition~\ref{prop:pre-Thomfund1} and Theorem~\ref{thm:plus-minus-map-flow}. We then prove Theorem~\ref{thm:plus-minus-map-flow} which is the analogue for the maps of $\T^{gl}$ of Theorem~\ref{thm:plus-minus-map} for the globular subdivisions. It is the third key fact of the paper. 

Section~\ref{section:final} uses all preceding results to prove that the image by the functor $\dcat:\ptop{\mathcal{M}} \to \dtop$ of every globular subdivision factors as a map of $\cell(\T^{cof})$ followed by a weak equivalence of flows. We obtain Corollary~\ref{final} which is the goal of the paper. The section is concluded with a discussion about the link between globular subdivision and underlying homotopy type of flow and ends with Corollary~\ref{cor:final2}.

Appendix~\ref{section:underlying-space} establishes Ken Brown's lemma for the underlying space functor from multipointed $d$-spaces to topological spaces in Proposition~\ref{prop:almost-kenbrown} despite the fact that it is not a left Quillen adjoint.

\subsection*{Acknowledgments}

I thank the referee for the reading and for the observations which have been incorporated into the final version.

\section{Categorical preliminaries}
\label{section:cat}

The knowledge of \cite{hocont}, in particular the exact notion of a dihomotopy of flows, is not required to understand this paper. On the other hand, the results expounded in this work rely on \cite{Moore1,Moore2,Moore3,GlobularNaturalSystem,leftproperflow}. Reminders are given throughout the paper. 

The initial (final resp.) object of a category, if it exists, is denoted by $\varnothing$ ($\mathbf{1}$ resp.). The identity map of an object $X$ is denoted by $\id_X$. The notation $S.X$, where $S$ is a set and $X$ is an object of a cocomplete category, is the sum $\coprod_{s\in S}X$. 

All facts of this section are well-known or are variants of well-known facts. The arguments are sketched for the ease of the reader (this also enables us to fix the notations).

\begin{prop} \label{prop:about-pushout-squares}
Let $\K$ be a cocomplete category. Consider the commutative diagram of objects of $\K$
\[\begin{tikzcd}[row sep=3em, column sep=3em]
	A \arrow[r]\arrow[d] \arrow[phantom,"\underline{C}",pos=0.5,rd]  & B \arrow[r] \arrow[d] \arrow[phantom,"\underline{D}",pos=0.5,rd]& C \arrow[d] \\
	D \arrow[r] & E \arrow[r] & F
\end{tikzcd}\]
If $\underline{C}$ and $\underline{D}$ are pushout squares, then the composite square $\underline{C}+\underline{D}$ is a pushout square.
\end{prop}

\bpf The proof is straightforward. \epf

\begin{cor} \label{cor:about-pushout-squares}
	Let $\K$ be a cocomplete category. Consider the commutative diagram of objects of $\K$
	\[\begin{tikzcd}[row sep=3em, column sep=3em]
		A \arrow[r]\arrow[d] \arrow[phantom,"\underline{C}",pos=0.5,rd]  & B \arrow[r] \arrow[d] \arrow[phantom,"\underline{D}",pos=0.5,rd]& C \arrow[d] \\
		D \arrow[r] & E \arrow[r] & F
	\end{tikzcd}\]
	If $\underline{C}$ and $\underline{C}+\underline{D}$ are pushout squares, then the commutative square $\underline{D}$ is a pushout square.
\end{cor}

\bpf Assume that $\underline{C}$ and $\underline{C}+\underline{D}$ are pushout squares. Replace the commutative square $\underline{D}$ by a pushout square $\underline{D'}$. Then by Proposition~\ref{prop:about-pushout-squares}, $\underline{C}+\underline{D'}$ is a pushout square. Hence $\underline{C}+\underline{D} \iso \underline{C}+\underline{D'}$, $\underline{C}+\underline{D}$ being a pushout square. We deduce that the commutative square $\underline{D}$ is a pushout square.
\epf

Let $\K$ be a cocomplete category. A \textit{transfinite tower} in $\K$ consists of a colimit-preserving functor $X:\lambda\to \K$ from a transfinite ordinal $\lambda$ viewed as a small category (thanks to its poset structure) to $\K$. Let $X_\lambda=\liminj X$. The map $X_0\to X_\lambda$ is called the \textit{transfinite composition} of the transfinite tower. A transfinite tower is an object of the category $\K^\lambda$ of functors from $\lambda$ to $\K$.

The notation $f\boxslash g$ means that $f$ satisfies the \textit{left lifting property} (LLP) with respect to $g$, or equivalently that $g$ satisfies the \textit{right lifting property} (RLP) with respect to $f$. For a class of maps $\C$, let $\inj(\C) = \{g \in \K, \forall f \in \C, f\boxslash g\}$, $\cof(\C)=\{f\mid \forall g\in \inj(\C), f\boxslash g\}$ and $\cell(\C)$ denotes the class of transfinite compositions of pushouts of elements of $\C$. In a locally presentable category, and if $\C$ is a set of maps, $\cof(\C)$ is the class of maps which are retracts of maps of $\cell(\C)$ by \cite[Corollary~2.1.15]{MR99h:55031}. 

Let $\K$ be a bicomplete category.  Let  $\D:I\to\K$ be a diagram over a Reedy category $(I,I_+,I_-)$. The latching category at $i\in I$ is denoted by $\de(I_+\ddownarrow i)$,  the latching object at $i\in I$ by $L_i\D := \liminj_{\de(I_+\ddownarrow i)} \D$, the matching category at $i\in I$ by $\de(i\ddownarrow I_-)$ and the matching object at $i\in I$ by $M_i\D = \limproj_{\de(i\ddownarrow I_-)} \D$. See e.g. \cite[Definition~15.2.3 and Definition~15.2.5]{ref_model2} or \cite[Definition~5.1.2]{MR99h:55031} for the general definitions of a latching/matching category/object.

A model category is a bicomplete category $\K$ equipped with a class of cofibrations $\C$, a class of fibrations $\F$ and a class of weak equivalences $\W$ such that: 1) $\W$ is closed under retract and satisfies \ttt, 2) the pairs $(\C,\W\cap \F)$ and $(\C\cap \W,\F)$ are functorial weak factorization systems. We refer to \cite[Chapter~1]{MR99h:55031} and to \cite[Chapter~7]{ref_model2} for the basic notions about general model categories. A cofibration is denoted by \begin{tikzcd}[cramped]\bullet\arrow[r,rightarrowtail]&\bullet\end{tikzcd}, a fibration by \begin{tikzcd}[cramped]\bullet\arrow[r,twoheadrightarrow]&\bullet\end{tikzcd} and a weak equivalence by \begin{tikzcd}[cramped]\bullet\arrow[r,"\simeq"]&\bullet\end{tikzcd}. 

A \textit{cellular} object $X$ of a cofibrantly generated model category is an object such that the canonical map $\varnothing\to X$ belongs to $\cell(I)$ where $I$ is the set of generating cofibrations. The maps of $\cell(I)$ are called \textit{cellular maps}. The transfinite sequence of pushouts is called a \textit{cellular decomposition} of $X$ and each pushout is called a \textit{cell}. In a combinatorial model category (i.e. a model category such that the underlying category is locally presentable \cite{MR2506258}), every map $X\to Y$ factors as a composite $X\to Z\to Y$ such that the left-hand map belongs to $\cell(I)$ and the right-hand map to $\inj(I)$ by \cite[Proposition~1.3]{MR1780498}.

\begin{prop} \label{prop:homotopy-colimit-tower}
	Let $\K$ be a model category. Let $F,G:\lambda\to \K$ be two transfinite towers of cofibrant objects such that all maps of the towers are cofibrations. Let $\mu:F\Rightarrow G$ be an objectwise weak equivalence. Then $\liminj \mu:\liminj F\to \liminj G$ is a weak equivalence.
\end{prop}

\bpf
By \cite[Theorem~5.1.3]{MR99h:55031}, the colimit functor is a left Quillen functor if the category of towers $\K^\lambda$ is equipped with its Reedy model structure and the two towers are Reedy cofibrant. Hence $\liminj \mu:\liminj F\to \liminj G$ is a weak equivalence.
\epf

By replacing in the statement of Proposition~\ref{prop:hocolim-tower} \textit{cocomplete category} by \textit{model category} and $\cell(I)$ by \textit{cofibration}, the same argument as the one of Proposition~\ref{prop:homotopy-colimit-tower} completes the proof. Since Proposition~\ref{prop:hocolim-tower} is slightly more general, we borrow the argument from \cite[Lemma~9.3.4]{cofibrationcat} which is initially written in the setting of cofibration categories.

\begin{prop} \label{prop:hocolim-tower} (variant of \cite[Corollary~5.1.5]{MR99h:55031} and \cite[Lemma~9.3.4]{cofibrationcat})
Let $\K$ be a cocomplete category. Let $I$ be a set of maps of $\K$. Let $\lambda$ be a transfinite ordinal. Consider two transfinite towers $A:\lambda \to \C$ and $B:\lambda \to \C$ and a natural transformation $f:A\Rightarrow B$. Assume that for each ordinal $\nu<\lambda$, the map $B_\nu \sqcup_{A_\nu} A_{\nu+1} \to B_{\nu+1}$ belongs to $\cell(I)$. Then the map 
\[
B_0 \sqcup_{A_0} A_\lambda \longrightarrow B_\lambda
\]
belongs to $\cell(I)$ as well.
\end{prop}

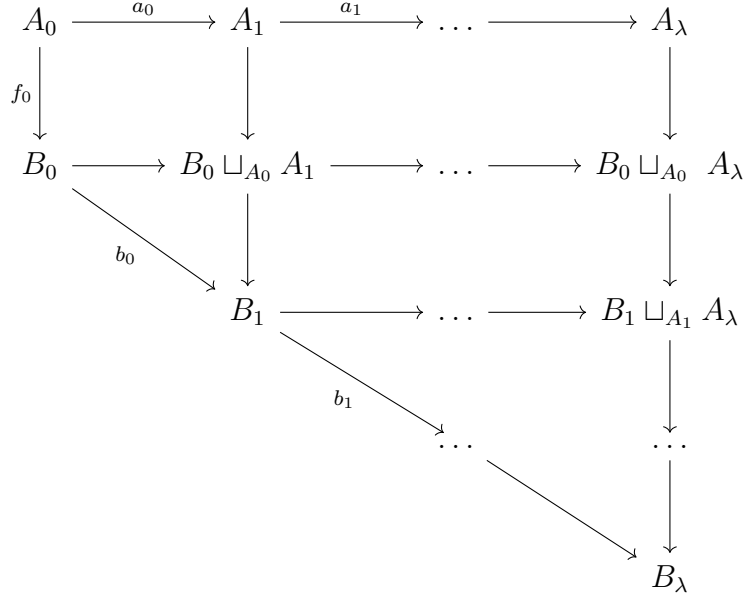
\begin{figure}
	\[
	\begin{tikzcd}[row sep=3em, column sep=3em]
		A_0 \arrow[r,"a_0"] \arrow[d,"f_0"'] & A_1 \arrow[r,"a_1"] \arrow[d] & \dots \arrow[r] & A_\lambda \arrow[d]\\
		B_0 \arrow[rd,"b_{0}"'] \arrow[r] &  B_0\sqcup_{A_0} A_1 \arrow[r]\arrow[d] & \dots \arrow[r] & B_0\sqcup_{A_0} A_\lambda \arrow[d] \\
		& B_1 \arrow[rd,"b_1"'] \arrow[r] & \dots \arrow[r] & B_1 \sqcup_{A_1} A_\lambda \arrow[d] \\
		&&\dots \arrow[rd] & \dots \arrow[d]\\
		&&& B_\lambda
	\end{tikzcd}
	\]
	\caption{Colimits of maps of $\cell(I)$}
	\label{fig:colim-cell-maps}
\end{figure}

\bpf
We consider the commutative diagram of $\K$ of Figure~\ref{fig:colim-cell-maps}. The map $B_0 \sqcup_{A_0} A_\lambda \to B_\lambda$ is the transfinite composition of the maps $B_{\nu} \sqcup_{A_{\nu}} A_\lambda \to B_{\nu+1} \sqcup_{A_{\nu+1}} A_\lambda$ with $\nu<\lambda$, and each map of the transfinite sequence is a pushout of $B_{\nu}\sqcup_{A_{\nu}} A_{\nu+1} \to B_{\nu+1}$ along the map $B_{\nu} \sqcup_{A_{\nu}} A_{\nu+1} \to B_{\nu} \sqcup_{A_{\nu}} A_\lambda$. The proof is complete, each map $B_{\nu}\sqcup_{A_{\nu}} A_{\nu+1} \to B_{\nu+1}$ belonging to $\cell(I)$ by hypothesis.
\epf

\begin{prop} \label{prop:useful-lemma}
	Let $\K$ be a model category. Consider a commutative diagram of solid arrows 
	\[
	\begin{tikzcd}[row sep=3em, column sep=5em]
		&&\\
		A \arrow[d,rightarrowtail] \arrow[dashed,r,"\ell"] \arrow[rr,bend left=30pt]  & B \arrow[d]\ \arrow[r,twoheadrightarrow,"\simeq"] & C \arrow[d] \\
		A' \arrow[dashed,r,"\ell'"] \arrow[rr,bend right=30pt] & B' \arrow[twoheadrightarrow,r,"\simeq"] & C'\\
		&&
	\end{tikzcd}
	\]
	such that the maps $B\to C$ and $B'\to C'$ are trivial fibrations and such that the map $A\to A'$ is a cofibration between cofibrant objects. Then there exist two maps $\ell:A\to B$ and $\ell':A'\to B'$ making the diagram commutative.
\end{prop}

\bpf
Consider the Reedy model category $\K^{1\to 2}$ of functors from the direct Reedy category $1\to 2$ to $\K$. It coincides with the projective model structure since the Reedy category $1\to 2$ is direct. In this model category of maps, the map $A\to A'$ is Reedy cofibrant, and the commutative square 
\[
\begin{tikzcd}[row sep=3em, column sep=5em]
	B \arrow[d] \arrow[twoheadrightarrow,"\simeq",r] & C \arrow[d] \\
	B' \arrow[r,twoheadrightarrow,"\simeq"] & C'
\end{tikzcd}
\]
is a trivial fibration from the map $B\to B'$ to the map $C\to C'$. The commutative diagram 
\[
\begin{tikzcd}[row sep=3em, column sep=10em]
	A \arrow[d,rightarrowtail] \arrow[r] &   C \arrow[d] \\
	A'  \arrow[r]&  C'
\end{tikzcd}
\]
is a map in $\K^{1\to 2}$ from the cofibrant object $A\to A'$ to the object $C\to C'$. Hence the existence of a lift $(\ell,\ell')$ from $A\to A'$ to $B\to B'$ making the whole diagram in the statement of the proposition commutative.
\epf

\section{Multipointed d-space and flow}
\label{section:space-flow}

The category $\top$ denotes either the category of \textit{$\Delta$-generated spaces} or of \textit{$\Delta$-Hausdorff $\Delta$-generated spaces} (cf. \cite[Section~2 and Appendix~B]{leftproperflow}). It is Cartesian closed by a result due to Dugger and Vogt recalled in \cite[Proposition~2.5]{mdtop} and locally presentable by \cite[Corollary~3.7]{FR}. The internal hom is denoted by $\ttop(-,-)$. The right adjoint of the inclusion from $\Delta$-generated spaces to general topological spaces is called the $\Delta$-kelleyfication. We will make use of the well-known three model structures of $\top$, namely the q-model structure, the h-model structure and the m-model structure: see the end of \cite[Appendix~B]{leftproperflow} for an overview and the bibliography of the latter paper for many other references (e.g. \cite{mixed-cole}).

Let $\gamma_1$ and $\gamma_2$ be two continuous maps from $[0,1]$ to some topological space such that $\gamma_1(1)=\gamma_2(0)$. The composite defined by 
\[
(\gamma_1 *_N \gamma_2)(t) = \begin{cases}
	\gamma_1(2t)& \hbox{ if }0\leq t\leq \frac{1}{2},\\
	\gamma_2(2t-1)& \hbox{ if }\frac{1}{2}\leq t\leq 1
\end{cases}
\]
is called the \textit{normalized composition}. 

Let $\mathcal{M}$ be the set of non decreasing surjective maps from $[0,1]$ to $[0,1]$ equipped with the $\Delta$-kelleyfication of the relative topology induced by the set inclusion $\mathcal{M} \subset \ttop([0,1],[0,1])$ (this set is denoted $\mathcal{M}(1,1)$ in \cite{Moore1,Moore3,GlobularNaturalSystem}).

\bd \cite[Definition~3.4]{Moore3} \label{def:multipointed-d-space} A \textit{multipointed $d$-space $X$} is a triple $(|X|,X^0,\P^{top}X)$ such that
\begin{itemize}
	\item The pair $(|X|,X^0)$ is a \textit{multipointed space}. The space $|X|$ is called the \textit{underlying space} of $X$ and the set $X^0$ the \textit{set of states} of $X$.
	\item The set $\P^{top}X$ is a set of continuous maps from $[0,1]$ to $|X|$ called the \textit{execution paths}, satisfying the following axioms:
	\begin{itemize}
		\item For any execution path $\gamma$, one has $\gamma(0),\gamma(1)\in X^0$.
		\item Let $\gamma$ be an execution path of $X$. Then any composite $\gamma\phi$ with $\phi\in \mathcal{M}$ is an execution path of $X$.
		\item Let $\gamma_1$ and $\gamma_2$ be two composable execution paths of $X$; then the normalized composition $\gamma_1 *_N \gamma_2$ is an execution path of $X$.
	\end{itemize}
\end{itemize}
A map $f:X\to Y$ of multipointed $d$-spaces is a map of multipointed spaces from $(|X|,X^0)$ to $(|Y|,Y^0)$ such that for any execution path $\gamma$ of $X$, the map $\P^{top}f:\gamma\mapsto f. \gamma$ is an execution path of $Y$. The category of multipointed $d$-spaces is denoted by $\ptop{\mathcal{M}}$. Let $\P_{\alpha,\beta}^{top} X = \{\gamma\in \P^{top}X\mid \gamma(0)=\alpha,\gamma(1)=\beta\}$. The set $\P^{top}_{\alpha,\beta} X$ is equipped with the $\Delta$-kelleyfication of the relative topology with respect to the inclusion $\P^{top}_{\alpha,\beta} X \subset \ttop([0,1],|X|)$. Thus a set map $Z\to \P^{top}_{\alpha,\beta} X$ where $Z$ is $\Delta$-generated, is continuous if and only if the associated map $Z\p [0,1] \to |X|$ is continuous. 
\ed

The category $\ptop{\mathcal{M}}$ is locally presentable, and in particular bicomplete, by \cite[Proposition~3.6]{Moore3}. The functor which forgets the execution paths $\Omega:X\mapsto (|X|,X^0)$ from $\ptop{\mathcal{M}}$ to the category $\mtop$ of multipointed topological spaces is topological in the sense of \cite[Section~21]{topologicalcat} by \cite[Theorem~3.17]{Moore3}. In particular, it creates limits and colimits thanks to the initial and final structures respectively. 

Every set $S$ can be viewed as a multipointed $d$-spaces $(S,S,\varnothing)$. The \textit{topological globe of a topological space $Z$}, which is denoted by $\globM(Z)$, is the multipointed $d$-space defined as follows
\begin{itemize}
	\item the underlying topological space is the quotient space \[\frac{\{{0},{1}\}\sqcup (Z\p[0,1])}{(z,0)=(z',0)={0},(z,1)=(z',1)={1}}\]
	\item the set of states is $\{{0},{1}\}$
	\item the set of execution paths is the set of continuous maps \[\{\delta_z\phi\mid \phi\in \mathcal{M},z\in  Z\}\]
	with $\delta_z(t) = (z,t)$.	It is equal to the underlying set of the space $Z \p \mathcal{M}$.
\end{itemize}
In particular, $\globM(\varnothing)$ is the multipointed $d$-space $\{{0},{1}\} = (\{{0},{1}\},\{{0},{1}\},\varnothing)$. 

\begin{nota}
	Let $n\geq 1$. Denote by $\mathbf{D}^n = \{(x_1,\dots,x_n)\in \mathbb{R}^n, x_1^2+\dots + x_n^2 \leq 1\}$ the $n$-dimensional disk, and by $\mathbf{S}^{n-1} = \{(x_1,\dots,x_n)\in \mathbb{R}^n, x_1^2+\dots + x_n^2 = 1\}$ the $(n-1)$-dimensional sphere. By convention, let $\mathbf{D}^{0}=\{0\}$ and $\mathbf{S}^{-1}=\varnothing$. 
\end{nota}

\begin{nota}
	Denote by $\vI^{top} = \globM(\mathbf{D}^0)$ the directed segment multipointed with the extremities.
\end{nota}

\bd \label{def:restriction-multidspace}
Let $X=(|X|,X^0,\P^{top}X)$ be a multipointed $d$-space. Let $A\subset X^0$. The \textit{restriction of $X$ to $A$} is the multipointed $d$-space 
\[
X\rest_{A} = (|X|,A,\{\gamma\in \P^{top} X\mid \gamma(0),\gamma(1)\in A\}).
\]
\ed

The \textit{q-model structure} of multipointed $d$-spaces \cite[Section~4]{Moore3} is the unique combinatorial model structure such that 
\[ I^{gl,top}\cup \{C:\varnothing \to \{0\},R:\{0,1\} \to \{0\}\}\]
with $I^{gl,top}=\{\globM(\mathbf{S}^{n-1})\subset \globM(\mathbf{D}^{n}) \mid n\geq 0\}$ is the set of generating cofibrations, the maps between globes being induced by the closed inclusions $\mathbf{S}^{n-1}\subset \mathbf{D}^{n}$, and such that 
\[
J^{gl,top}=\{\globM(\mathbf{D}^{n})\subset \globM(\mathbf{D}^{n+1}) \mid n\geq 0\}
\]
is the set of generating trivial cofibrations, the maps between globes being induced by the closed inclusions $(x_1,\dots,x_n)\mapsto (x_1,\dots,x_n,0)$. The weak equivalences are the maps of multipointed $d$-spaces $f:X\to Y$  inducing a bijection $f^0:X^0\iso Y^0$ and a weak homotopy equivalence $\P^{{top}} f:\P^{{top}}_{\alpha,\beta} X \to \P^{{top}}_{f(\alpha),f(\beta)} Y$ for all $(\alpha,\beta)\in X^0\p X^0$ and the fibrations are the maps of multipointed $d$-spaces $f:X\to Y$  inducing a q-fibration $\P_{\alpha,\beta}^{{top}} f:\P_{\alpha,\beta}^{{top}} X \to \P_{f(\alpha),(\beta)}^{{top}} Y$ of topological spaces for all $(\alpha,\beta)\in X^0\p X^0$. All multipointed $d$-spaces are q-fibrant.

\bd \cite[Definition~4.11]{model3} \label{def:flow}
A \textit{flow} $X$ is a small enriched semicategory. Its set of objects (preferably called \textit{states}) is denoted by $X^0$ and the space of morphisms (preferably called \textit{execution paths}) from $\alpha$ to $\beta$ is denoted by $\P_{\alpha,\beta}X$ (e.g. \cite[Definition~10.1]{Moore1}). For any $x\in \P_{\alpha,\beta}X$, let $s(x)=\alpha$ and $t(x)=\beta$. The category is denoted by $\dtop$. Let 
\[
\P X = \coprod_{(\alpha,\beta)\in X^0\p X^0} \P_{\alpha,\beta} X.
\]
\ed

The category $\dtop$ is locally presentable. Every set can be viewed as a flow with an empty path space. This give rise to a functor from sets to flows which is limit-preserving and colimit-preserving. More generally, any poset can be viewed as a flow, with a unique execution path from $u$ to $v$ if and only if $u<v$. This gives rise to a functor from the category of posets together with the strictly increasing maps to flows.

\begin{nota} \label{nota:glob}
	For any topological space $Z$, the flow $\glob(Z)$ is the flow having two states $0$ and $1$ and such that the only nonempty space of execution paths, when $Z$ is nonempty, is $\P_{0,1}\glob(Z)=Z$. It is called \textit{the globe of $Z$}. Let $\vI = \glob(\mathbf{D}^0)$.
\end{nota}

The \textit{q-model structure} of flows \cite[Theorem~7.6]{QHMmodel} is the unique combinatorial model structure such that 
\[I^{gl} \cup \{C:\varnothing \to \{0\},R:\{0,1\} \to \{0\}\}\]
with $I^{gl}=\{\glob(\mathbf{S}^{n-1})\subset \glob(\mathbf{D}^{n}) \mid n\geq 0\}$ is the set of generating cofibrations, the maps between globes being induced by the closed inclusions $\mathbf{S}^{n-1}\subset \mathbf{D}^{n}$, and such that 
\[
J^{gl}=\{\glob(\mathbf{D}^{n})\subset \glob(\mathbf{D}^{n+1}) \mid n\geq 0\}
\]
is the set of generating trivial cofibrations, the maps between globes being induced by the closed inclusions $(x_1,\dots,x_n)\mapsto (x_1,\dots,x_n,0)$. The weak equivalences are the maps of flows $f:X\to Y$  inducing a bijection $f^0:X^0\iso Y^0$ and a weak homotopy equivalence  $\P f:\P_{\alpha,\beta} X \to \P_{f(\alpha),f(\beta)} Y$ for all $(\alpha,\beta)\in X^0\p X^0$ and the fibrations are the maps of flows $f:X\to Y$  inducing a q-fibration $\P_{\alpha,\beta} f:\P_{\alpha,\beta} X \to \P_{f(\alpha),(\beta)} Y$ of topological spaces for all $(\alpha,\beta)\in X^0\p X^0$. All flows are q-fibrant.

\bth \label{thm:complement-flow} (\cite[Theorem~5.7]{leftproperflow}) Let $f:X\to Y$ be a (trivial resp.) q-cofibration between q-cofibrant flows. Then the continuous map $\P f:\P X \to \P Y$ is a (trivial resp.) q-cofibration of spaces. In particular, the path space functor $\P :\dtop\to \top$ preserves q-cofibrancy. 
\eth

Both in the q-model categories of multipointed $d$-spaces and flows, for any cellular map $f:X\to Y$, the cells $C:\varnothing \subset \{0\}$ and $R:\{0,1\}\to \{0\}$ can be regrouped at the beginning of the cellular decomposition of $f$ as follows (since we only need this fact for the q-model structure of multipointed $d$-spaces in this paper, the proof for the category of flows, which is actually similar, is omitted)~\footnote{I used the existence of this factorization many time in my previous papers without proving it.}.

\begin{prop} \label{prop:regroupCR}
	Let $f:X\to Y$ be a cellular map of multipointed $d$-spaces. The map $f$ factors as a composite \[f:X\longrightarrow Z\longrightarrow Y\] such that the map $X\to Z$ belongs to $\cell(\{C,R\})$ and such that the map $Z\to Y$ is bijective on states (or equivalently belongs to $\inj(\{C,R\})$) and belongs to  $\cell(I^{gl,top})$, and finally such that the transfinite sequence of cells belonging to $I^{gl,top}$ in the cellular decomposition of $f$ is equal to the transfinite sequence of cells of the cellular decomposition of $Z\to Y$. Moreover, the factorization is unique up to isomorphism.
\end{prop}

\bpf By \cite[Proposition~1.3]{MR1780498}, the map $f$ factors as a composite $f:X\to Z\to Y$ such that the map $g:X\to Z$ belongs to $\cell(\{C,R\})$ and such that the map $h:Z\to Y$ belongs to $\inj(\{C,R\})$. This implies that the map $h:Z\to Y$ is bijective on states. By hypothesis, there is a transfinite tower $T:\lambda\to \ptop{\mathcal{M}}$ of pushouts of maps of $I^{gl,top}\cup\{C,R\}$ such that the transfinite composition $T_0\to T_\lambda$ is the map $f$. We are going to build a transfinite tower $T':\lambda\to \ptop{\mathcal{M}}$ of pushouts of maps of $I^{gl,top}$ such that the transfinite composition is the map $h:Z\to Y$ as follows. Let \[g_0=g:T_0=X\longrightarrow T'_0=Z.\] Let $\nu<\lambda$. Assume that there is a commutative diagram of multipointed $d$-spaces 
\[
\begin{tikzcd}[row sep=3em,column sep=3em]
	X\arrow[d,"g"] \arrow[r] & T_\nu \arrow[d,"g_\nu"]\arrow[r] & T_\lambda=Y \arrow[d,equal]\\
	Z  \arrow[r]\arrow[rr,bend right=30pt,"h"] & T'_\nu \arrow[r] & T'_\lambda=Y
\end{tikzcd}
\]
such that the bottom horizontal maps are bijective on states. There are three mutually exclusive cases. 
\begin{enumerate}[leftmargin=*]
	\item The map $T_\nu \to T_{\nu+1}$ is a pushout of the map $C:\varnothing\to \{0\}$. Then there is a unique map $g_{\nu+1}$ making the following diagram commutative
	\[
	\begin{tikzcd}[row sep=3em,column sep=3em]
		X\arrow[d,"g"] \arrow[r] & T_\nu \arrow[d,"g_\nu"]\arrow[r] & T_{\nu+1} \arrow[dashed,d,"g_{\nu+1}"]\arrow[r] & T_\lambda=Y \arrow[d,equal]\\
		Z  \arrow[r] \arrow[rrr,bend right=30pt,"h"] & T'_\nu \arrow[r,equal] & T'_{\nu+1} \arrow[r] & T'_\lambda=Y
	\end{tikzcd}
	\]
	Indeed, one has $|T_{\nu+1}|=|T_\nu|\sqcup \{\alpha\}$. The additional state $\alpha\in T_{\nu+1}^0\backslash T_\nu^0$ is taken to some state $\alpha'$ of $T_\lambda=T'_\lambda=Y$. Since the bottom horizontal maps are bijective on states, $\alpha'\in Z^0$ and the image by the map $Z\to T'_\nu$ yields $g_{\nu+1}(\alpha)$.
	\item The map $T_\nu \to T_{\nu+1}$ is a nontrivial pushout of the map $R:\{0,1\}\to \{0\}$. Then there is a unique map $g_{\nu+1}$ making the following diagram commutative
	\[
	\begin{tikzcd}[row sep=3em,column sep=3em]
		X\arrow[d,"g"] \arrow[r] & T_\nu \arrow[d,"g_\nu"]\arrow[r] & T_{\nu+1} \arrow[dashed,d,"g_{\nu+1}"]\arrow[r] & T_\lambda=Y \arrow[d,equal]\\
		Z  \arrow[r]\arrow[rrr,bend right=30pt,"h"] & T'_\nu \arrow[r,equal] & T'_{\nu+1} \arrow[r] & T'_\lambda=Y
	\end{tikzcd}
	\]
	Indeed, there exist two states $\alpha_1$ and $\alpha_2$ of $T_\nu$ such that $|T_{\nu+1}|=|T_\nu|/(\alpha_1=\alpha_2)$. The two states $\alpha_1$ and $\alpha_2$ of $T_\nu$ are taken to some state $\alpha'$ of $T_\lambda=T'_\lambda=Y$. Since the bottom horizontal maps are bijective on states, $\alpha'\in Z^0$ and the image by the map $Z\to T'_\nu$ gives $g_{\nu+1}(\alpha)$. Consequently, the continuous map $|g_\nu|:|T_\nu|\to |T'_\nu|$ factors uniquely as a composite $|T_\nu|\to |T_{\nu+1}| \to |T'_\nu| = |T'_{\nu+1}|$. The right-hand continuous map gives rise to a map of multipointed $d$-spaces $g_{\nu+1}$ since the execution paths of $T_{\nu+1}$ are generated by the normalized compositions of execution paths of $T_{\nu}$ generated by the identification $\alpha_1=\alpha_2$.
	\item The map $T_\nu \to T_{\nu+1}$ is a pushout of a map of the form $\globM(\mathbf{S}^{n-1})\to \globM(\mathbf{D}^n)$ for some $n\geq 0$. Then the map $g_{\nu+1}$ is defined by a pushout diagram as follows
	\[
	\begin{tikzcd}[row sep=3em,column sep=3em]
		X\arrow[d,"g"] \arrow[r] & T_\nu \arrow[d,"g_\nu"]\arrow[r] & T_{\nu+1} \arrow[d,"g_{\nu+1}"]\arrow[r] & T_\lambda=Y \arrow[d,equal]\\
		Z  \arrow[r]\arrow[rrr,bend right=30pt,"h"] & T'_\nu \arrow[r,] & \cocartesian T'_{\nu+1} \arrow[r] & T'_\lambda=Y
	\end{tikzcd}
	\]
\end{enumerate}
For a limit ordinal $\nu\leq \lambda$, $g_\nu$ is defined as \[g_\nu=\liminj_{\mu<\nu}g_\mu.\] We obtain the commutative diagram of multipointed $d$-spaces 
\[
\begin{tikzcd}[row sep=3em,column sep=3em]
	X\arrow[d,"g"] \arrow[r] & T_\lambda \arrow[d,"g_\lambda"]\arrow[r,equal]  & T_\lambda=Y \arrow[d,equal]\\
	Z  \arrow[r]\arrow[rr,bend right=30pt,"h"] & T'_\lambda \arrow[r,equal] & T'_\lambda=Y
\end{tikzcd}
\]
Hence the proof of the existence is complete. Two factorizations give rise to a commutative diagram of solid arrows of multipointed $d$-spaces
\[
\begin{tikzcd}[column sep=3em, row sep=3em]
	X\arrow[d,"g"'] \arrow[dd,bend right=50pt,"f"']\arrow[rd,"g'"]\arrow[r,equal] & X \arrow[d,"g'"]\arrow[dd,bend left=50pt,"f"]\\
	Z \arrow[dashed,"\ell_{Z,Z'}",r] \arrow[rd,"h"]\arrow[d,"h"']& Z'\arrow[d,"h'"]\\
	Y \arrow[r,equal] & Y
\end{tikzcd}
\]
Since $g\in\cell(\{C,R\})$ and $h'\in \inj(\{C,R\})$, the lift $\ell_{Z,Z'}$ exists. Take another lift $\ell'_{Z,Z'}$. Then $h'\ell_{Z,Z'}=h=h'\ell'_{Z,Z'}$. Since $h'\in \cell(I^{gl,top})$, the map $|h'|:|Z'|\to |Y|$ is a q-cofibration of spaces by Proposition~\ref{prop:double-suspension}. Therefore it is one-to-one. We deduce that $\ell_{Z,Z'}=\ell'_{Z,Z'}$. We obtain that $\id_Z=\ell_{Z,Z}$ for all $Z$ realizing the factorization and therefore that $\id_Z=\ell_{Z,Z'}\ell_{Z',Z}$ and $\id_{Z'}=\ell_{Z',Z}\ell_{Z,Z'}$. We have proved the uniqueness of the factorization up to isomorphism.
\epf

\section{The categorization functor}
\label{section:dcat}

There is a unique functor $\dcat:\ptop{\mathcal{M}}\longrightarrow \dtop$ from the category of multipointed $d$-spaces to the category of flows, called the \textit{categorization functor}, taking a multipointed $d$-space $X$ to the unique flow $\dcat(X)$ such that $\dcat(X)^0=X^0$ and such that $\P_{\alpha,\beta}\dcat(X)$ is the quotient of the space of execution paths $\P_{\alpha,\beta}^{top}X$ by the equivalence relation generated by the reparametrization by $\mathcal{M}$, the composition of $\dcat(X)$ being induced by the normalized composition. One has $\dcat(\vI^{top})=\vI$. For any topological space $Z$, there is the natural isomorphism of flows $\dcat(\globM(Z))\iso \glob(Z)$.

The functor $\dcat:\ptop{\mathcal{M}}\longrightarrow \dtop$ is neither a left adjoint nor a right adjoint \cite[Proposition~7.3]{mdtop} and \cite[Theorem~8.8]{Moore2}: there is a slight difference in the definition of multipointed $d$-space used in the latter theorems which does not affect the proofs.

\begin{prop} \label{prop:rappel-all}
There is an isomorphism of functors $\dcat \iso \lmoore\moore^{top}$ such that $\lmoore$ is some left adjoint and $\moore^{top}$ some right adjoint such that for all transfinite towers $X:\lambda \to \ptop{\mathcal{M}}$ of maps of $I^{gl,top}$ such that $X_0=(S,S,\varnothing)$ for some set $S$ equipped with the discrete topology, the natural map \[\liminj_{\nu<\lambda} \moore^{top}(X_\nu)\longrightarrow \moore^{top}(\liminj_{\nu<\lambda} X_\nu)\] is an isomorphism in the target category of the functor $\moore^{top}$.
\end{prop}  

\bpf There is an isomorphism $\dcat \iso \lmoore\moore^{top}$ by \cite[Notation~19 and Theorem~15]{Moore3}. The second statement is \cite[Theorem~5]{Moore3}.
\epf

\bth \label{thm:rappel-all}
The functor $\dcat:\ptop{\mathcal{M}}\longrightarrow \dtop$ satisfies the following properties:
\begin{enumerate}
	\item Its left derived functor $\mathbf{L}\dcat$ in the sense of \cite{HomotopicalCategory} induces an equivalence of categories between the homotopy categories of the q-model structures of multipointed $d$-spaces and flows.
	\item The functor $\dcat$ takes weak equivalences between q-cofibrant multipointed $d$-spaces to weak equivalences between q-cofibrant flows.
	\item Let $X:\lambda\to\ptop{\mathcal{M}}$ be a transfinite tower of pushouts of maps of $I^{gl,top}\cup\{C,R\}$ between cellular multipointed $d$-spaces. Then the canonical map \[\liminj_{\nu<\lambda} \dcat(X_\nu)\longrightarrow \dcat(\liminj_{\nu<\lambda} X_\nu)\] is an isomorphism of flows.
	\item Let $X:\lambda\to\ptop{\mathcal{M}}$ be a transfinite tower of pushouts of maps of $\cell(I^{gl,top}\cup\{C,R\})$ between cellular multipointed $d$-spaces. Then the canonical map \[\liminj_{\nu<\lambda} \dcat(X_\nu)\longrightarrow \dcat(\liminj_{\nu<\lambda} X_\nu)\] is an isomorphism of flows.
\end{enumerate} 
\eth

\bpf $(1)$ and $(2)$ is \cite[Theorem~15]{Moore3}. Let us treat $(3)$ now. Consider a cellular decomposition $X':\lambda'\to \ptop{\mathcal{M}}$ of the cellular multipointed $d$-space $X_0$ such that $X'_0=(S',S',\varnothing)$ for some set $S'$ equipped with the discrete topology and such that all cells belong to $I^{gl,top}$. Such a cellular decomposition always exists by regrouping the pushouts of the maps $C:\varnothing\to\{0\}$ and $R:\{0,1\}\to \{0\}$ of the tower $X'$ at the beginning of the tower $X'$ using Proposition~\ref{prop:regroupCR}. We obtain a transfinite tower denoted by $X'+X:\lambda'+\lambda\to \ptop{\mathcal{M}}$, where $\lambda'+\lambda$ is the ordinal sum, such that $(X'+X)_0=(S',S',\varnothing)$ and such that the colimit is equal to the colimit of $X:\lambda\to\ptop{\mathcal{M}}$. Using Proposition~\ref{prop:regroupCR}, we regroup the pushouts of the maps $C:\varnothing\to\{0\}$ and $R:\{0,1\}\to \{0\}$ of the tower $X$ at the beginning of the tower $X'+X$. We obtain another transfinite tower $X'':\lambda''\to\ptop{\mathcal{M}}$ with $\lambda''\leq \lambda'+\lambda$ having the same colimit as $X$ such that $X''_0=(S'',S'',\varnothing)$ for another set $S''$ equipped with the discrete topology and such that each map of the tower $X''$ belongs to $I^{gl,top}$. Using Proposition~\ref{prop:rappel-all}, the natural map \[\liminj_{\nu<\lambda''} \moore^{top}(X''_\nu)\longrightarrow \moore^{top}(\liminj_{\nu<\lambda''} X''_\nu)\] is an isomorphism in the target category of the functor $\moore^{top}$. We obtain the isomorphisms of flows 
\begin{align*}
	\liminj_{\nu<\lambda} \dcat(X_\nu) & \iso \liminj_{\nu<\lambda''} \dcat(X''_\nu) \\
	& \iso \lmoore\bigg(\liminj_{\nu<\lambda''} \moore^{top}(X''_\nu)\bigg)\\
	&\iso \lmoore\bigg(\moore^{top}(\liminj_{\nu<\lambda''} X''_\nu)\bigg)\\
	& \iso \dcat(\liminj_{\nu<\lambda''} X''_\nu)\\
	& \iso \dcat(\liminj_{\nu<\lambda} X_\nu),
\end{align*}
the first isomorphism since the towers $\dcat.X$ and $\dcat.X''$ have same colimits by construction of $X''$, the second isomorphism since $\dcat \iso \lmoore\moore^{top}$ and since $\lmoore$ is a left adjoint, the third isomorphism by Proposition~\ref{prop:rappel-all}, the fourth isomorphism since  $\dcat \iso \lmoore\moore^{top}$ and finally the fifth isomorphism since the towers $X$ and $X''$ have same colimits by construction of $X''$. The proof of $(3)$ is complete. The implication $(3)\Rightarrow (4)$ is clear.
\epf

Corollary~\ref{cor:almost-accessible} is not used in the paper. However it is worth being mentioned.

\begin{cor} \label{cor:almost-accessible}
	Let $X:\lambda\to\ptop{\mathcal{M}}$ be a transfinite tower of q-cofibrations between q-cofibrant multipointed $d$-spaces. Then the canonical map \[\liminj_{\nu<\lambda} \dcat(X_\nu)\longrightarrow \dcat(\liminj_{\nu<\lambda} X_\nu)\] is an isomorphism of flows.
\end{cor}

\bpf
The functor category $\ptop{\mathcal{M}}^\lambda$ can be equipped with its Reedy model structure. The transfinite tower $X:\lambda\to\ptop{\mathcal{M}}$ is Reedy q-cofibrant by the definition of the Reedy model structure. Since the transfinite ordinal $\lambda$ equipped with its poset structure gives rise to a direct small category, the Reedy model structure coincides with the projective model structure. Thus the generating Reedy q-cofibrations are the maps of transfinite towers $\lambda(\nu,-).U\to \lambda(\nu,-).V$ with $\nu<\lambda$ such that $U\to V$ is a generating q-cofibration of the q-model category of multipointed $d$-spaces. Adding the cell $\lambda(\nu,-).U\to \lambda(\nu,-).V$ to a cellular object of $\ptop{\mathcal{M}}^\lambda$ consists of adding the cell $U\to V$ to each object of the tower starting from the ordinal $\nu<\lambda$, colimits being calculated objectwise in $\ptop{\mathcal{M}}^\lambda$ and since the category $\lambda$ is a poset. This implies, colimits commuting between one another, that a cellular object of $\ptop{\mathcal{M}}^\lambda$ is always a transfinite tower, and not just a tower. Factor the Reedy q-cofibration $\varnothing \to X$ as a composite $\begin{tikzcd}[cramped,column sep=small]\varnothing\arrow[r,rightarrowtail]&X^{cell} \arrow[r,twoheadrightarrow,"\simeq"]& X\end{tikzcd}$ of a cellular map of $\ptop{\mathcal{M}}^\lambda$ followed by a trivial Reedy q-fibration. The tower $X^{cell}$ is a transfinite tower as in Theorem~\ref{thm:rappel-all}~$(4)$. Therefore the canonical map \[\liminj_{\nu<\lambda} \dcat(X^{cell}_\nu)\longrightarrow \dcat(\liminj_{\nu<\lambda} X^{cell}_\nu)\] is an isomorphism of flows. The lift $\ell$ 
in the diagram of towers 
\[
\begin{tikzcd}[row sep=3em,column sep=3em]
	& X^{cell} \arrow[d,twoheadrightarrow,"\simeq"]\\
	X\arrow[r,equal] \arrow[ru,dashed,"\ell"]& X
\end{tikzcd}
\]
exists since the tower $X$ is Reedy q-cofibrant. Hence $X$ is a retract of $X^{cell}$. The proof is complete since a retract of an isomorphism is an isomorphism.
\epf

\begin{prop}\label{prop:colim} Consider a pushout diagram of multipointed $d$-spaces of the form 
	\[
	\begin{tikzcd}[row sep=3em, column sep=4em]
		U \arrow[r] \arrow[d,"i"'] & A \arrow[d] \\
		V \arrow[r] & \cocartesian B
	\end{tikzcd}
	\]
	where $i:U\to V \in I^{gl,top}\cup\cell(\{C\})$ with $A$ cellular and $n\geq 0$. Then there is a pushout diagram of flows 
	\[
	\begin{tikzcd}[row sep=3em, column sep=3em]
		\dcat(U) \arrow[r] \arrow[d,"I^{gl}\cup\cell(\{C\})\ni"'] & \dcat(A) \arrow[d] \\
		\dcat(V) \arrow[r] & \cocartesian \dcat(B).
	\end{tikzcd}
	\]
\end{prop}

\bpf The case $i\in I^{gl}$ is treated in \cite[Proposition~7.1]{GlobularNaturalSystem}. The case $i\in\cell(\{C\})$ is different: it comes from the fact that $B=(S,S,\varnothing)\sqcup A$ for some set $S$. This implies that $\dcat(B)=S\sqcup \dcat(A)$ by definition of the functor $\dcat:\ptop{\mathcal{M}}\to \dtop$ where the set $S$ is viewed as a discrete poset, and hence as a flow.
\epf

\bth \label{thm:colim}
Consider a pushout diagram of cellular multipointed $d$-spaces 
\[
\begin{tikzcd}[row sep=3em, column sep=3em]
	A \arrow[r] \arrow[d,rightarrowtail,"\cell(I^{gl,top}\cup\{C\})\ni"'] \arrow[rd,phantom,"\underline{C}"]& E \arrow[d,rightarrowtail] \\
	B \arrow[r] & \cocartesian F
\end{tikzcd}
\]
such that the vertical maps belong to $\cell(I^{gl,top}\cup\{C\})$. Then there is the pushout diagram of cellular flows 
\[
\begin{tikzcd}[row sep=3em, column sep=3em]
	\dcat(A) \arrow[r] \arrow[d,rightarrowtail,"\cell(I^{gl}\cup\{C\})\ni"'] & \dcat(E) \arrow[d,rightarrowtail] \\
	\dcat(B) \arrow[r] & \cocartesian \dcat(F)
\end{tikzcd}
\]
and the vertical maps belong to $\cell(I^{gl}\cup\{C\})$.
\eth

\bpf
Let $i\in I^{gl}\cup\{C\}$. Assume first that there is a pushout diagram of cellular multipointed $d$-spaces
\[
\begin{tikzcd}[row sep=3em, column sep=3em]
	U \arrow[r] \arrow[d,rightarrowtail,"i"'] & A \arrow[r] \arrow[d,rightarrowtail] & E \arrow[d,rightarrowtail] \\
	V \arrow[r]  & \cocartesian B \arrow[r] & \cocartesian F
\end{tikzcd}
\]
for some $n\geq 0$. By Proposition~\ref{prop:about-pushout-squares}, there is the pushout square of cellular multipointed $d$-spaces 
\[
\begin{tikzcd}[row sep=3em, column sep=3em]
	U \arrow[r] \arrow[d,rightarrowtail,"i"'] &  E \arrow[d,rightarrowtail] \\
	V \arrow[r]  &  \cocartesian F
\end{tikzcd}
\]
Using Proposition~\ref{prop:colim}, we obtain the pushout squares 
\[
\begin{tikzcd}[row sep=3em, column sep=3em]
	\dcat(U) \arrow[r] \arrow[d,rightarrowtail,"i"'] &  \dcat(A) \arrow[d,rightarrowtail] & \dcat(U) \arrow[r] \arrow[d,rightarrowtail,"\dcat(i)"'] &  \dcat(E) \arrow[d,rightarrowtail] \\
	\dcat(V) \arrow[r]  &  \cocartesian \dcat(B) & \dcat(V) \arrow[r]  &  \cocartesian \dcat(F)
\end{tikzcd}
\]
Using Corollary~\ref{cor:about-pushout-squares}, we deduce the pushout diagram of cellular flows 
\[
\begin{tikzcd}[row sep=3em, column sep=3em]
	\dcat(A) \arrow[r] \arrow[d,rightarrowtail] & \dcat(E) \arrow[d,rightarrowtail] \\
	\dcat(B) \arrow[r] & \cocartesian \dcat(F)
\end{tikzcd}
\]
We then deduce the theorem when the vertical maps of $\underline{C}$ are a pushout of a finite composition of maps of $I^{gl,top}\cup\{C\}$ using Proposition~\ref{prop:about-pushout-squares}. The transfinite case is a consequence of Theorem~\ref{thm:rappel-all}~$(4)$. All involved flows are cellular by Proposition~\ref{prop:colim} and Theorem~\ref{thm:rappel-all}~$(4)$.
\epf

We do not know whether Proposition~\ref{prop:colim} and Theorem~\ref{thm:colim} hold by considering the cofibration $R:\{0,1\}\to \{0\}$. When a pushout of $R$ in the category of multipointed $d$-spaces is nontrivial, new execution paths are created. Using the theory of Moore flows developed in \cite{Moore1,Moore2,Moore3} and in particular their structure of biclosed semimonoidal category, it is easy to describe the underlying set of the space of execution paths. It is not clear what its topology is. Appendix~\ref{section:underlying-space} shows that the case of this cofibration must be treated with caution.

\section{Enriched small category and directed space}
\label{section:enriched-id}

To have simpler statements for Theorem~\ref{thm:plus-minus-map-flow}, we introduce the category of topologically enriched small categories.

\begin{nota} \label{def:II}
	The category of (topologically) enriched small categories is denoted by $\cat_\top$. The forgetful functor $\cat_\top\subset \dtop$ has a left adjoint denoted by $\II:\dtop\to\cat_\top$. 
\end{nota}

The execution path of a flow $X$ from $\alpha$ to $\beta$ is still denoted by $\P_{\alpha,\beta}X$, whereas the space of morphisms from $\alpha$ to $\beta$ in $\II(X)$ is denoted by $\II(X)(\alpha,\beta)$. One has the homeomorphisms 
\[
\II(X)(\alpha,\beta) = \begin{cases}
	\P_{\alpha,\beta}X & \hbox{ if }\alpha\neq \beta\\
	\{\id_\alpha\}\sqcup \P_{\alpha,\beta}X& \hbox{ if }\alpha= \beta.
\end{cases}
\]

We also need to introduce Grandis' notion of directed space to be able to use some results from \cite{GlobularNaturalSystem}.

\begin{nota}
	Denote by $\mathcal{I}$ the set of non-decreasing continuous maps from $[0,1]$ to $[0,1]$. Note that an element of $\mathcal{I}$ can be a constant map (this set is denoted $\mathcal{I}(1)$ in \cite{GlobularNaturalSystem}).
\end{nota}

\bd \cite[Definition~1.1]{mg} \cite[Definition~4.1]{DAT_book} \label{def:directed_space}
A \textit{directed space} is a pair $X=(|X|,\vec{P}(X))$ consisting of a topological space $|X|$ and a set $\vec{P}(X)$ of continuous paths from $[0,1]$ to $|X|$ satisfying the following axioms:
\begin{itemize}
	\item $\vec{P}(X)$ contains all constant paths;
	\item $\vec{P}(X)$ is closed under normalized composition;
	\item $\vec{P}(X)$ is closed under reparametrization by an element of $\mathcal{I}$.
\end{itemize}
The space $|X|$ is called the \textit{underlying topological space} or the \textit{state space}. The elements of $\vec{P}(X)$ are called \textit{directed paths}. A morphism of directed spaces is a continuous map between the underlying topological spaces which takes a directed path of the source to a directed path of the target. The category of directed spaces is denoted by $\ptop{}$. Write $\vec{P}(X)(u,v)$ for the space of directed paths of $X$ from $u$ to $v$ equipped with the $\Delta$-kelleyfication of the compact-open topology.
\ed

The \textit{category of traces} of a directed space $X$, denoted by $\vec{T}(X)$, has for objects the points of $X$ and the set of maps $\vec{T}(X)(a,b)$ from $a\in X$ to $b\in X$ is the set of traces $\tr{\gamma}$ of directed paths $\gamma$ going from $a$ to $b$, i.e. the set of directed paths from $a$ to $b$ up to reparametrization by a map of $\mathcal{M}$. The composition of traces, denoted by $*$, is induced by the normalized composition of directed paths, i.e. $\tr{\gamma} * \tr{\gamma'} = \tr{\gamma *_N\gamma'}$. It is strictly associative. 

By \cite[Proposition~3.7 and Theorem~3.8]{GlobularNaturalSystem}, the mapping $\discont:Y=(|Y|,\vec{P}(Y)) \mapsto (|Y|,|Y|,\vec{P}(Y))$ induces a full and faithful functor $\discont:\ptop{}\to\ptop{\mathcal{M}}$ which is a right adjoint. 

\begin{nota} \label{notation:cont}
	Denote by $\cont:\ptop{\mathcal{M}} \to \ptop{}$ the left adjoint.
\end{nota}

There is the equality $\cont(\discont(X))=X$ for all directed spaces $X$. By \cite[Proposition~3.6]{GlobularNaturalSystem}, the left adjoint $\cont:\ptop{\mathcal{M}} \to \ptop{}$ is defined as follows. The underlying space of $\cont(X)$ is $|X|$ and the set of directed spaces $\vec{P}(X)$ consists of all constant paths and all Moore compositions of the form $[(\gamma_1\phi_1\mu_{\ell_1}) * \dots * (\gamma_n\phi_n\mu_{\ell_n})$ such that $\ell_1+\dots + \ell_n = 1$ where $\gamma_1,\dots,\gamma_n$ are execution paths of $X$ and $\phi_i\in \mathcal{I}$ for $i=1,\dots,n$, and where $\mu_\ell:[0,\ell]\to [0,1]$ is defined by $\mu_\ell(t)=t/\ell$ with $\ell>0$.

\bth \label{thm:continuous} 
Let $X$ be a cellular multipointed $d$-space for the q-model category of $\ptop{\mathcal{M}}$. A continuous map from $[0,1]$ to $|X|$ is a directed path of $X$, i.e. of $\cont(X)$, if and only if it is of the form $\gamma\phi$ where $\gamma$ is an execution path of $X$ or a constant path and where $\phi\in \mathcal{I}$. 
\eth

\bpf
It is  \cite[Theorem~4.9]{GlobularNaturalSystem}.
\epf

\begin{nota}
	For a multipointed $d$-space $X$, let $\vec{T}(X)=\vec{T}(\cont(X))$. For all $\alpha,\beta\in |X|$, let $\vec{P}(X)(\alpha,\beta)=\vec{P}(\cont(X))(\alpha,\beta)$.
\end{nota}

\begin{prop} \label{prop:trace-pathspace}
Let $X$ be a q-cofibrant multipointed $d$-space $X$. Let $\alpha,\beta\in X^0$. Then there are the homeomorphisms 
\[
\vec{P}(X)(\alpha,\beta) \iso \begin{cases}
	\{\alpha\} \sqcup \P^{top}_{\alpha,\beta}X & \hbox{ if }\alpha=\beta\\
	\P^{top}_{\alpha,\beta}X & \hbox{ if }\alpha\neq \beta
\end{cases}
\]
where $\alpha$ denotes the constant path $\alpha$. There is also the homeomorphism
\[
\vec{T}(X)(\alpha,\beta) \iso \II(\dcat(X))(\alpha,\beta)
\]
for all $\alpha,\beta\in X^0$. Moreover the spaces $\vec{P}(X)(\alpha,\beta)$ are m-cofibrant and the spaces $\vec{T}(X)(\alpha,\beta)$ are q-cofibrant.
\end{prop}

\bpf
Every q-cofibrant multipointed $d$-space is a retract of a cellular one. Thus one can suppose that $X$ is cellular. The first part is then a consequence of Theorem~\ref{thm:continuous}. The second part is the consequence of the definitions of a trace and of the functor $\dcat:\ptop{\mathcal{M}}\to\dtop$. The spaces $\P^{top}_{\alpha,\beta}X$ are m-cofibrant by \cite[Theorem~16]{Moore3} and the spaces $\P_{\alpha,\beta}\dcat(X)$ are q-cofibrant by Theorem~\ref{thm:complement-flow}, $\dcat(X)$ being a q-cofibrant flow by Theorem~\ref{thm:rappel-all}~$(2)$. Hence the proof is complete.
\epf

\begin{figure}
	\[
\begin{tikzpicture}[black,scale=2.5,pn/.style={circle,inner sep=0pt,minimum width=4pt,fill=dark-red}]
	\draw (0,0) node[pn] {} node[black,below left] {$(0\,,0)$};
	\draw (1,0) node[pn] {} node[black,below right] {$(1\,,0)$};
	\draw (0,1) node[pn] {} node[black,above left] {$(0\,,1)$};
	\draw (1,1) node[pn] {} node[black,above right] {$(1\,,1)$};
	\draw[->] [very thick] (0.0,0.0) -- (0.5,0.5);
	\draw[->] [very thick] (0.5,0.5) -- (1.0,0.0);
	\draw[->] (0,0) -- (1,1);
	\draw[->] (0,1) -- (1,0); 
\end{tikzpicture}
\]
\caption{The fake crossing}
\label{fig:fakecrossing}
\end{figure}
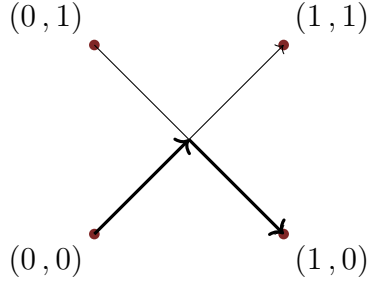

Consider the fake crossing $X$ depicted in Figure~\ref{fig:fakecrossing} (it already appears in \cite[Figure~2]{GlobularNaturalSystem}). The underlying space $|X|$ is the topological space  $|X| = \{(u,u)\mid u\in [0,1]\} \cup \{(u,1-u)\mid u\in [0,1]\}$. Let $X^0 = \{(0,0),(0,1),(1,0),(1,1)\}$. Let $\P^{top}_{(0,0),(1,1)}X = \{t\mapsto (\phi(t),\phi(t))\mid \phi\in \mathcal{M}\}$, $\P^{top}_{(0,1),(1,0)}X = \{t\mapsto (\phi(t),1-\phi(t))\mid \phi\in \mathcal{M}\}$ and $\P^{top}_{\alpha,\beta}X=\varnothing$ otherwise. Then $\vec{P}(X)((0,0),(1,0)) \iso \mathcal{M}$ and $\P_{(0,0),(1,0)}^{top}X=\varnothing$. The multipointed $d$-space $X$ is not q-cofibrant because a q-cofibrant replacement of $X$ is the disjoint sum $\vI^{top}\sqcup \vI^{top}$ of two directed segments multipointed with the extremities.
	
For every precubical set $K$, consider its realization $|K|$ as a multipointed $d$-space: the set of states is $K_0$, the underlying space is the geometric realization of $K$ and the set of execution path from $\alpha$ to $\beta$ is the set of all nonconstant directed paths in the geometric realization of $K$ from $\alpha$ to $\beta$ (\cite[Definition~3.12]{GlobularNaturalSystem}). This yields a functor from precubical sets to multipointed $d$-spaces. It satisfies trivially the homeomorphisms 
\[
\vec{P}(|K|)(\alpha,\beta) \iso \begin{cases}
	\{\alpha\} \sqcup \P^{top}_{\alpha,\beta}|K| & \hbox{ if }\alpha=\beta\\
	\P^{top}_{\alpha,\beta}|K| & \hbox{ if }\alpha\neq \beta
\end{cases}
\]
and 
\[
\vec{T}(|K|)(\alpha,\beta) \iso \II(\dcat(|K|))(\alpha,\beta)
\]
for all $\alpha,\beta\in K_0$.

This implies that the q-cofibrancy hypothesis is a sufficient but not necessary condition to obtain the conclusion of Proposition~\ref{prop:trace-pathspace}.

\bth \label{thm:carac-weakeq-trace}
	Consider a map of multipointed $d$-spaces $f:X\to Y$ between two q-cofibrant multipointed $d$-spaces. The following statements are equivalent: 
	\begin{enumerate}
		\item $f$ is a weak equivalence of the q-model structure.
		\item $f$ induces a bijection between the states and for all $(\alpha,\beta)\in X^0\p X^0$ a weak homotopy equivalence $\vec{P}(X)(\alpha,\beta) \simeq \vec{P}(Y)(f(\alpha),f(\beta))$.
		\item $f$ induces a bijection between the states and for all $(\alpha,\beta)\in X^0\p X^0$ a weak homotopy equivalence $\vec{T}(X)(\alpha,\beta) \simeq \vec{T}(Y)(f(\alpha),f(\beta))$.
		\item $f$ induces a bijection between the states and for all $(\alpha,\beta)\in X^0\p X^0$ a homotopy equivalence $\vec{P}(X)(\alpha,\beta) \simeq \vec{P}(Y)(f(\alpha),f(\beta))$.
		\item $f$ induces a bijection between the states and for all $(\alpha,\beta)\in X^0\p X^0$ a homotopy equivalence $\vec{T}(X)(\alpha,\beta) \simeq \vec{T}(Y)(f(\alpha),f(\beta))$.
	\end{enumerate} 
\eth

\bpf
The equivalence $(1)\Leftrightarrow (2)$ is a consequence of the definition of a weak equivalence and of Proposition~\ref{prop:trace-pathspace}. By Proposition~\ref{prop:trace-pathspace} and \cite[Theorem~16]{Moore3}, the quotient maps $\vec{P}(X)(\alpha,\beta)\to \vec{T}(X)(\alpha,\beta)$ and $\vec{P}(Y)(f(\alpha),f(\beta))\to \vec{T}(Y)(f(\alpha),f(\beta))$ are homotopy equivalences. For all $\alpha,\beta\in X^0$, there is the commutative diagram of spaces 
\[
\begin{tikzcd}[column sep=3em,row sep=3em]
	\vec{P}(X)(\alpha,\beta) \arrow[r]\arrow[d,"\simeq"'] & \vec{P}(Y)(f(\alpha),f(\beta))\arrow[d,"\simeq"]\\
	\vec{T}(X)(\alpha,\beta) \arrow[r] & \vec{T}(Y)(f(\alpha),f(\beta))
\end{tikzcd}
\]
The equivalence $(2)\Leftrightarrow (3)$ is then a consequence of the \ttt. Since $X$ and $Y$ are q-cofibrant, the topological spaces $\vec{P}(X)(\alpha,\beta)$ and $\vec{P}(Y)(f(\alpha),f(\beta))$ are m-cofibrant by Proposition~\ref{prop:trace-pathspace}. Moreover, the topological spaces $\vec{T}(X)(\alpha,\beta)$ and $\vec{T}(Y)(f(\alpha),f(\beta))$ are q-cofibrant by Proposition~\ref{prop:trace-pathspace} as well. We obtain the equivalences $(2) \Leftrightarrow (4)$ and  $(3) \Leftrightarrow (5)$ thanks to \cite[Corollary~3.4]{mixed-cole}.
\epf

\section{Cellular multipointed d-space}
\label{section:cellular}

All cellular multipointed $d$-spaces for the q-model structure of $\ptop{\mathcal{M}}$ can be reached from $\varnothing$ without using the cofibration $R:\{0,1\}\to \{0\}$ and by regrouping the pushouts of $C:\varnothing\to \{0\}$ at the very beginning using Proposition~\ref{prop:regroupCR}. Thus, for the sequel, a cellular decomposition of a cellular multipointed $d$-space of the q-model category $\ptop{\mathcal{M}}$ consists of a colimit-preserving functor $X:\lambda \longrightarrow \ptop{\mathcal{M}}$ from a transfinite ordinal $\lambda$ to the category of multipointed $d$-spaces such that
\begin{itemize}
	\item The multipointed $d$-space $X_0$ is a set, in other terms $X_0=(X^0,X^0,\varnothing)$ for some set $X^0$.
	\item For all $\nu<\lambda$, there is a pushout diagram of multipointed $d$-spaces 
	\[\begin{tikzcd}[row sep=3em, column sep=3em]
		\globM(\mathbf{S}^{n_\nu-1}) \arrow[d,rightarrowtail] \arrow[r,"g_\nu"] & X_\nu \arrow[d,rightarrowtail] \\
		\globM(\mathbf{D}^{n_\nu}) \arrow[r,"\widehat{g_\nu}"] & \cocartesian X_{\nu+1}
	\end{tikzcd}\]
	with $n_\nu \geq 0$. 
\end{itemize}
The underlying topological space $|X_\lambda|$ is Hausdorff by \cite[Proposition~4.4]{Moore3}. For all $\nu\leq \lambda$, there is the equality $X_\nu^0=X^0$. Denote by \[c_\nu = |\globM(\mathbf{D}^{n_\nu})|\backslash |\globM(\mathbf{S}^{n_\nu-1})|\] the $\nu$-th cell of $X_\lambda$. It is called a \textit{globular cell}. Like in the usual setting of CW-complexes, $\widehat{g_\nu}$ induces a homeomorphism from $c_\nu$ to $\widehat{g_\nu}(c_\nu)$ equipped with the relative topology. The map $\widehat{g_{\nu}}: \globM(\mathbf{D}^{n_\nu})\to X_\lambda$ is called the \textit{attaching map} of the globular cell $c_\nu$. The state $\widehat{g_\nu}(0)\in X^0$ ($\widehat{g_\nu}(1)\in X^0$ resp.)  is called the \textit{initial (final resp.) state} of $c_\nu$ and is denoted by $c_\mu^-$ ($c_\mu^+$ resp.). The integer $n_\nu+1$ is called the \textit{dimension} of the globular cell $c_\nu$. It is denoted by $\dim c_\nu$. The states of $X^0$ are also called the \textit{globular cells of dimension $0$}. By convention, a state of $X^0$ viewed as a globular cell of dimension $0$ is equal to its initial state and to its final state. Thus, for $\alpha\in X^0$, one has $\alpha=\alpha^+=\alpha^-$. The set of globular cells of $X_\lambda$ is denoted by $\mathcal{C}(X_\lambda)$. The set of globular cells of dimension $n\geq 0$ of $X_\lambda$ is denoted by $\mathcal{C}_n(X_\lambda)$. In particular, $\mathcal{C}_0(X_\lambda)=X^0$. The closure $\widehat{g_\nu}(c_\nu)$ of $c_\nu$ in $|X_\lambda|$ is denoted by $\widehat{c_{\nu}}$.

\begin{nota}
	Denote by $\cell_f$ the class of finite compositions of pushouts of the inclusions $\{\mathbf{S}^{n-1}\subset \mathbf{D}^n\mid n\geq 0\}$. Let 
	\begin{align*}
		&\mathbf{S}^{n-1}_0 = \{(x_1,\dots,x_n,0)\mid \sum_i x_i^2=1\}, \\
		& \mathbf{D}^n_0 = \{(x_1,\dots,x_n,0)\mid \sum_i x_i^2\leq 1\}.
	\end{align*}
	There are the homeomorphisms $\mathbf{S}^{n-1}_0 \iso \mathbf{S}^{n-1}$ and $\mathbf{D}^n_0 \iso \mathbf{D}^n$ for all $n\geq 0$.
\end{nota}

\begin{prop} \label{prop:multi-well-pointed-2}
Let $n\geq 1$. Let $F$ be a finite subset of the interior of $\mathbf{D}^n$. Then the inclusion $\mathbf{S}^{n-1} \cup F\subset \mathbf{D}^n$ belongs to $\cell_f$.
\end{prop}

\bpf
For $n=1$, write $F=\{u_1<\dots <u_p\}$ with $p\geq 1$. One has \[\mathbf{D}^1 = [-1,1] =\bigcup_{1\leq k\leq p+1} [u_{k-1},u_k]\] with $u_0=-1$ and $u_{p+1}=1$. This implies that $\mathbf{S}^0 \cup F\subset \mathbf{D}^1$ belongs to $\cell_f$. We prove now by induction on $n\geq 1$ 
\[
\mathcal{E}(n):\forall u\in \mathbf{D}^n\backslash \mathbf{S}^{n-1}, \mathbf{S}^{n-1} \cup \{u\} \subset \mathbf{D}^n \in \cell_f.
\]
The case $\mathcal{E}(1)$ is treated above. Assume $\mathcal{E}(n)$ for $n\geq 1$. We want to prove $\mathcal{E}(n+1)$. Using a homeomorphism, we can suppose that $u$ is the center of $\mathbf{D}^{n+1}$. Consider the commutative diagram of topological spaces 
\[
\begin{tikzcd}[row sep=3em, column sep=3em]
	\mathbf{S}^{n-1}_0\cup\{u\} \arrow[r,"\subset"] \arrow[d,"\subset"]& \mathbf{S}^n \cup \{u\} \arrow[d,"\subset"]\\
	\mathbf{D}^n_0 \arrow[r,"\subset"] & \cocartesian \mathbf{D}^n_0 \cup \mathbf{S}^n \arrow[d,"\subset"]\\
	& \mathbf{D}^{n+1}
\end{tikzcd}
\]
Since the square is a pushout, $\mathcal{E}(n)$ implies that the inclusion $\mathbf{S}^n \cup \{u\} \subset \mathbf{D}^n_0 \cup \mathbf{S}^n$ belongs to $\cell_f$. The inclusion $\mathbf{D}^n_0 \cup \mathbf{S}^n\subset \mathbf{D}^{n+1}$ belongs to $\cell_f$ as well since $\mathbf{D}^{n+1}$ is obtained from $\mathbf{D}^n_0 \cup \mathbf{S}^n$ by using two pushouts along the inclusion $\mathbf{S}^n\subset \mathbf{D}^{n+1}$. We have proved $\mathcal{E}(n+1)$. We proceed now by induction on $p\geq 0$ to prove 
\[
\mathcal{E'}(p) : \forall q\leq p, \forall n\geq 2, \mathbf{S}^{n-1} \cup \{u_1,\dots ,u_q\} \subset \mathbf{D}^n \in \cell_f.
\]
There is nothing to prove for $\mathcal{E'}(0)$ and $\mathcal{E'}(1)$ is already proved above. Assume $\mathcal{E'}(p)$ for $p\geq 1$. We want to prove $\mathcal{E'}(p+1)$. Using a homeomorphism, one can suppose that $u_{p+1}$ is the center of $\mathbf{D}^n$ and that $\mathbf{D}^{n-1}_0 \cap \{u_1,\dots,u_p\} = \varnothing$. Consider the commutative diagram of topological spaces 
\[
\begin{tikzcd}[row sep=3em, column sep=3em]
	\mathbf{S}^{n-2}_0\cup\{u\} \arrow[r,"\subset"] \arrow[d,"\subset"]& \mathbf{S}^{n-1} \cup \{u_1,\dots,u_{p+1}\} \arrow[d,"\subset"]\\
	\mathbf{D}^{n-1}_0 \arrow[r,"\subset"] & \cocartesian \mathbf{D}^{n-1}_0 \cup \mathbf{S}^{n-1} \cup\{u_1,\dots,u_p\}\arrow[d,"\subset"]\\
	& \mathbf{D}^{n}
\end{tikzcd}
\]
Since the square is a pushout, $\mathcal{E}(n-1)$ implies that the inclusion $\mathbf{S}^{n-1} \cup \{u_1,\dots,u_{p+1}\} \subset \mathbf{D}^{n-1}_0 \cup \mathbf{S}^{n-1} \cup\{u_1,\dots,u_p\}$ belongs to $\cell_f$. The inclusion $\mathbf{D}^{n-1}_0 \cup \mathbf{S}^{n-1} \cup\{u_1,\dots,u_p\}\subset \mathbf{D}^{n}$ belongs to $\cell_f$ as well since $\mathbf{D}^{n}$ is obtained from $\mathbf{D}^{n-1}_0 \cup \mathbf{S}^{n-1} \cup\{u_1,\dots,u_p\}$ by using two pushouts along two inclusions of the form $\mathbf{S}^{n-1} \cup F\subset \mathbf{D}^{n}$ with $F\subset \mathbf{D}^{n}\backslash \mathbf{S}^{n-1}$ and $F$ of cardinal lower than $p$. We have proved $\mathcal{E'}(p+1)$.
\epf

\begin{prop} \label{prop:subd-gl-globe}
Let $n\geq 0$. Consider a finite set \[F\subset |\globM(\mathbf{D}^n)| \backslash |\globM(\mathbf{S}^{n-1})|.\] Then the following data assemble into a multipointed $d$-space denote by $\globM(\mathbf{D}^n)_F$:
\begin{itemize}
	\item The set of states is $\{0,1\} \cup F$.
	\item The underlying space is $|\globM(\mathbf{D}^n)|$.
	\item For all $\alpha\neq \beta \in \{0,1\} \cup F$, $\P^{top}_{\alpha,\beta} \globM(\mathbf{D}^n)_F = \vec{P}(\cont(\globM(\mathbf{D}^n)))(\alpha,\beta)$. 
	\item For all $\alpha\in \{0,1\} \cup F$, $\P^{top}_{\alpha,\alpha} \globM(\mathbf{D}^n)_F=\varnothing$.  
\end{itemize}
In particular, there is the isomorphism of multipointed $d$-spaces \[\globM(\mathbf{D}^n)\iso \globM(\mathbf{D}^n)_\varnothing.\] 
\end{prop}

\bpf
The composition of two execution paths is an execution path and the set of execution paths is closed under reparametrization by $\mathcal{M}$.
\epf

\begin{prop} \label{prop:poset}
	Let $n\geq 0$. Consider a finite set \[F\subset |\globM(\mathbf{D}^n)| \backslash |\globM(\mathbf{S}^{n-1})|.\] The finite set $\{0,1\} \cup F$ can be equipped with a poset structure as follows: $u<v$ if and only if $\P_{u,v}^{top}\globM(\mathbf{D}^n)_F$ is nonempty (which implies contractible in this particular case). Moreover, the identity of $\globM(\mathbf{D}^n)_F^0=\{0,1\} \cup F$ induces a trivial q-fibration of flows \[\begin{tikzcd} \glob(\mathbf{D}^n)_F \arrow[r,twoheadrightarrow,"\simeq"]&(\{0,1\}\cup F,\leq).\end{tikzcd}\]
\end{prop}

\bpf
The map of flows $\glob(\mathbf{D}^n)_F\to (\{0,1\}\cup F,\leq)$ is a q-fibration since all spaces of execution paths of the poset $(\{0,1\}\cup F,\leq)$ are discrete. It is a weak equivalence since $\P_{u,v}^{top}\globM(\mathbf{D}^n)_F$  is either for $u\geq v$  empty or homeomorphic to $\mathcal{M}$ for $u<v$ which is contractible.
\epf

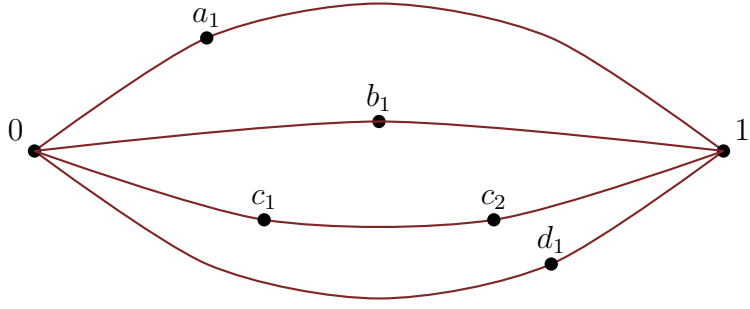
\begin{figure}
	\begin{tikzpicture}[scale=1.3,pn/.style={circle,inner sep=0pt,minimum width=5pt,fill=black}]
		\def\L{7}
		\draw[dark-red][thick] plot[smooth] coordinates {(0,0) (\L/4,1.15) (\L/2,1.5) (3*\L/4,1.15) (\L,0)};
		\draw[dark-red][thick] plot[smooth] coordinates {(0,0) (\L/4,-1.15) (\L/2,-1.5) (3*\L/4,-1.15) (\L,0)};
		\draw (0,0) node[pn] {} node[black,above left]{$0$};
		\draw (\L,0) node[pn] {} node[black,above right]{$1$};
		\draw (\L/2,0.3) node[pn] {} node[black,above]{$b_1$};
		\draw (\L/4,1.15) node[pn] {} node[black,above]{$a_1$};
		\draw (3*\L/4,-1.15) node[pn] {} node[black,above]{$d_1$};
		\draw[dark-red][thick] plot[smooth] coordinates {(0,0) (\L/2,0.3) (\L,0)};
		\def\courbe {plot[smooth] coordinates {(0,0) (\L/3,-0.7) (2*\L/3,-0.7) (\L,0)}}
		\draw[dark-red][thick] \courbe;
		\draw (\L/3,-0.7) node[pn] {} node[black,above]{$c_1$};
		\draw (2*\L/3,-0.7) node[pn] {} node[black,above]{$c_2$};
	\end{tikzpicture}
	\caption{The globular poset $(\{0,1\}\cup \{a_1,b_1,c_1,c_2,d_1\},\leq)$ with $0<a_1<1$, $0<b_1<1$, $0<c_1<c_2<1$ and $0<d_1<1$}
	\label{fig:glposet}
\end{figure}

The poset $(\{0,1\} \cup F,\leq)$ of Proposition~\ref{prop:poset} looks as follows: a smallest element $0$, a biggest element $1$ and a finite number of finite strictly increasing chains of the form $0<a_1<\dots <a_p<1$ with $p\geq 0$. One of these posets is depicted in Figure~\ref{fig:glposet}. Such a poset is called a \textit{globular poset}.

\begin{nota}
	Let $\ell<\ell'$ be two real numbers. The multipointed $d$-space $\vec{[\ell,\ell']}$ is defined as follows: the underlying space is the segment $[\ell,\ell']$, the set of states is $\{\ell,\ell'\}$ and the set of execution paths is the set of nondecreasing surjective maps from $[0,1]$ to $[\ell,\ell']$. For all topological spaces $Z$, the unique map $Z\to \{0\}$ induces a map of multipointed $d$-spaces $\pi:\globM(Z)\to \globM(\{0\}) \iso \vec{[0,1]}$. 
\end{nota}

\begin{prop} \label{prop:raf-cell}
Let $n\geq 0$. Consider a finite set \[F\subset |\globM(\mathbf{D}^n)| \backslash |\globM(\mathbf{S}^{n-1})|.\] Then the map of multipointed $d$-spaces $\globM(\mathbf{S}^{n-1}) \subset \globM(\mathbf{D}^n)_F$ is a finite composition of pushouts of maps of $I^{gl,top}\cup\{C\}$. 
\end{prop}

\bpf
The case $n=0$ is trivial. Assume that $n\geq 1$. We consider the pushout diagram of multipointed $d$-spaces (the poset structure is defined in Proposition~\ref{prop:poset})
\[
\begin{tikzcd}[row sep=3em, column sep=3em]
	\{0,1\}\sqcup \displaystyle\coprod_{\substack{(u,v)\in (\{0,1\} \cup F)^2 \\u<v,]u,v[=\varnothing,(u,v)\neq (0,1)}} \{u,v\} \arrow[d,"\subset"] \arrow[r,"\substack{u\mapsto u\\v \mapsto v}"] & \{0,1\} \cup F \arrow[d]\\
	\globM(\mathbf{S}^{n-1}) \sqcup \displaystyle\coprod_{\substack{(u,v)\in (\{0,1\} \cup F)^2 \\u<v,]u,v[=\varnothing,(u,v)\neq (0,1)}} \vec{[\pi(u),\pi(v)]} \arrow[r] & \cocartesian X_F
\end{tikzcd}
\]
Since the functor $X\mapsto X^0$ from multipointed $d$-spaces to sets is colimit-preserving, we have $X_F^0 = \{0,1\} \cup F$. Note that the two sums in the pushout diagram above are taken on all pairs $(u,v)$ such that $u<v$, $]u,v[=\varnothing$ and $(u,v)\neq (0,1)$. The latter condition implies that the sums are empty in the case $F=\varnothing$. This implies that $X_\varnothing = \globM(\mathbf{S}^{n-1})$. Consider the equivalence relation on $F$ induced by $u\sim v$ if and only if $u<v$. Since $u<v$ implies that there is a directed path from $u$ to $v$, every element of the quotient set $F/\!\sim$ corresponds to a unique element of the interior of $\mathbf{D}^n$. Thus, there is a canonical inclusion $F/\!\sim \subset \mathbf{D}^n\backslash \mathbf{S}^{n-1}$. By induction on the cardinal of $F$, we can verify that the multipointed $d$-space $X_F\rest_{\{0,1\}}$ (see Definition~\ref{def:restriction-multidspace}) is isomorphic to $\globM(\mathbf{S}^{n-1} \sqcup F/\!\sim)$. There is a pushout diagram of multipointed $d$-spaces 
\[
\begin{tikzcd}[row sep=3em, column sep=3em]
	X_F\rest_{\{0,1\}} \arrow[r] \arrow[d,"\subset"'] & X_F \arrow[d] \\
	\globM(\mathbf{D}^{n}) \arrow[r] & \cocartesian \globM(\mathbf{D}^n)_F
\end{tikzcd}
\]
We can now conclude the proof. The inclusion $\globM(\mathbf{S}^{n-1}) \subset \globM(\mathbf{D}^n)_F$ is the composite of the three inclusions 
\begin{align*}
	&\globM(\mathbf{S}^{n-1}) \subset\globM(\mathbf{S}^{n-1}) \sqcup \displaystyle\coprod_{\substack{(u,v)\in (\{0,1\} \cup F)^2 \\u<v,]u,v[=\varnothing,(u,v)\neq (0,1)}} \vec{[\pi(u),\pi(v)]}\\
	& \globM(\mathbf{S}^{n-1}) \sqcup \displaystyle\coprod_{\substack{(u,v)\in (\{0,1\} \cup F)^2 \\u<v,]u,v[=\varnothing,(u,v)\neq (0,1)}} \vec{[\pi(u),\pi(v)]} \subset X_F \\
	& X_F \subset \globM(\mathbf{D}^{n})_F
\end{align*}
Thanks to Proposition~\ref{prop:multi-well-pointed-2}, we deduce that the inclusion $\globM(\mathbf{S}^{n-1}) \subset \globM(\mathbf{D}^n)_F$ is a finite composition of pushouts of the maps $R:\{0,1\}\to \{0\}$, $C:\varnothing\to\{0\}$ and $\globM(\mathbf{S}^{k-1}) \subset \globM(\mathbf{D}^{k})$ for $k\geq 0$. Since the inclusion $\globM(\mathbf{S}^{n-1}) \subset \globM(\mathbf{D}^n)_F$ is one-to-one on states, the map $R:\{0,1\}\to \{0\}$ can be removed. And the proof is complete.
\epf

\begin{cor} \label{cor:glob-subd}
	Let $n\geq 0$. Consider a finite set \[F\subset |\globM(\mathbf{D}^n)| \backslash |\globM(\mathbf{S}^{n-1})|.\] The multipointed $d$-space $\globM(\mathbf{D}^n)_F$ is cellular.
\end{cor}

\begin{nota}
	Let $\glob(\mathbf{D}^n)_F=\dcat(\globM(\mathbf{D}^n)_F)$. We have $\glob(\mathbf{D}^n)_\varnothing = \glob(\mathbf{D}^n)$ since $\glob(\mathbf{D}^n)=\dcat(\globM(\mathbf{D}^n))$. 
\end{nota}

\begin{cor} \label{cor:glob-subd-dcat}
	Let $n\geq 0$. Consider a finite set \[F\subset |\globM(\mathbf{D}^n)| \backslash |\globM(\mathbf{S}^{n-1})|.\] The map of flows $\glob(\mathbf{S}^{n-1}) \subset \glob(\mathbf{D}^n)_F$ is a finite composition of pushouts of maps of $I^{gl}\cup \{C\}$. In particular, it is a q-cofibration.
\end{cor}

\bpf
It is a consequence of Theorem~\ref{thm:colim}.
\epf

\section{Globular subdivision}
\label{section:gl-sbd}

\bd \cite[Definition~4.10]{diCW} \cite[Definition~9.1]{GlobularNaturalSystem} \label{def:globular-sbd}
A map of multipointed $d$-spaces $f:X\to Y$ is a \textit{globular subdivision} if both $X$ and $Y$ are cellular and if $f$ induces a homeomorphism between the underlying topological spaces of $X$ and $Y$. We say that $Y$ is a \textit{globular subdivision} of $X$ when there exists such a map. This situation is denoted by \[\begin{tikzcd}[cramped]f:X\arrow[r,"\sbd"]&Y.\end{tikzcd}\]
\ed

A very simple example of globular subdivision can be obtained by removing states from the distinguished set of states of a cellular multipointed $d$-space. Let $X=(|X|,X^0,\P^{top}X)$ be a cellular multipointed $d$-space. Let $A\subset X^0$. The identity of $|X|$ induces a globular subdivision $X\rest_{A}\to X$ as soon as $X\rest_{A}$ is cellular (see Definition~\ref{def:restriction-multidspace}).

We recall the important fact:

\bth \label{thm:T-homotopy} \cite[Theorem~9.3]{GlobularNaturalSystem}
Let $f:X\to Y$ be a map of multipointed $d$-spaces between two cellular multipointed $d$-spaces $X$ and $Y$. The following two conditions are equivalent: 
\begin{enumerate}
	\item The map $f:X\to Y$ is a globular subdivision.
	\item The map of directed spaces $\cont(f):\cont(X) \to \cont(Y)$ is an isomorphism.
\end{enumerate}
\eth

By adapting the proof of \cite[Theorem~4.12]{diCW}, we can easily verify that the categorical localization of the full subcategory of cellular multipointed $d$-spaces by the globular subdivisions is locally small: the cardinal of the set of maps in the categorical localization from $X$ to $Y$ is bounded by $\max(2^{\aleph_0},\# |X|,\# |Y|)$ where $\#S$ means the cardinal of $S$.

\begin{prop} \label{prop:sbd-path}
Let \begin{tikzcd}[cramped]f:X\arrow[r,"\sbd"]&Y\end{tikzcd} be a globular subdivision. For all $\alpha,\beta\in X^0$, there is the homeomorphisms $\P^{top}_{\alpha,\beta}X \iso \P^{top}_{f(\alpha),f(\beta)}Y$ and $\vec{T}(X)(\alpha,\beta) \iso \vec{T}(Y)(f(\alpha),f(\beta))$. 
\end{prop}

\bpf
By Theorem~\ref{thm:T-homotopy}, there is the isomorphism of directed spaces $\cont(f):\cont(X)\iso \cont(Y)$. We obtain the homeomorphisms $\vec{T}(X)(\alpha,\beta) \iso \vec{T}(Y)(f(\alpha),f(\beta))$ for all $\alpha,\beta\in X^0$. The homeomorphisms $\P^{top}_{\alpha,\beta}X \iso \P^{top}_{f(\alpha),f(\beta)}Y$ for all $\alpha,\beta\in X^0$ are a consequence of Proposition~\ref{prop:trace-pathspace}. 
\epf

\bd \label{def:plus-minus-map}
Let \begin{tikzcd}[cramped]f:X\arrow[r,"\sbd"]&Y\end{tikzcd} be a globular subdivision. Choose a cellular decomposition of $X$. Every point $\alpha\in Y^0$ is in a unique globular cell $c_\alpha$ of $X$. Let $\alpha^-=c_\alpha^-$ and $\alpha^+=c_\alpha^+$. Note that when $\alpha\in X^0$, then $\alpha^-=\alpha^+=\alpha$. We obtain two set maps $\alpha\mapsto \alpha^-$ and $\alpha\mapsto \alpha^+$ from $Y^0$ to $X^0$.
\ed

Proposition~\ref{prop:independence-choice-gl-subd} proves that the set maps $\alpha\mapsto \alpha^-$ and $\alpha\mapsto \alpha^+$ of Definition~\ref{def:plus-minus-map} from $Y^0$ to $X^0$ do not depend of the choice of the cellular decomposition of $X$.

\begin{prop} \label{prop:independence-choice-gl-subd}
Let \begin{tikzcd}[cramped]f:X\arrow[r,"\sbd"]&Y\end{tikzcd} be a globular subdivision. Consider two cellular decompositions $\mathcal{C}_0(X)$ and $\mathcal{C}_1(X)$ of $X$ and the four associated set maps $\alpha\mapsto \alpha_0^-$, $\alpha\mapsto \alpha_1^-$, $\alpha\mapsto \alpha_0^+$ and $\alpha\mapsto \alpha_1^+$ from $Y^0$ to $X^0$. Then for all $\alpha\in Y^0$, one has $\alpha_0^-=\alpha_1^-$ and $\alpha_0^+=\alpha_1^+$.
\end{prop}

\bpf We just have to consider the case $\alpha\in Y^0\backslash X^0$. The point $\alpha\in |X|=|Y|$ belongs to a unique globular cell $c_i$ of $\mathcal{C}_i(X)$ for $i=0,1$ with $\dim(c_0)\geq 1$ and $\dim(c_1)\geq 1$. There exists a unique execution path $\gamma_i$ of $X$ up to reparametrization from $c_i^-$ to $c_i^+$ for $i=0,1$ with $\gamma_i(]0,1[) \cap X^0 = \varnothing$. Thus $\alpha_0^-=c_0^-=c_1^-=\alpha_1^-$ and $\alpha_0^+=c_0^+=c_1^+=\alpha_1^+$.
\epf

\bd \label{def:connectionmaps}
Let \begin{tikzcd}[cramped]f:X\arrow[r,"\sbd"]&Y\end{tikzcd} be a globular subdivision. Then the maps $\alpha\mapsto \alpha^-$ and $\alpha\mapsto \alpha^+$ from $Y^0$ to $X^0$, which depend only on $f$ by Proposition~\ref{prop:independence-choice-gl-subd}, are called the \textit{connection maps} of the globular subdivision $f$.
\ed

Consider the identity $\id_X$ of a cellular multipointed $d$-space $X$. It is a globular subdivision. The connection maps of $\id_X$ are the identity of $X^0$. The notations $(-)^-,(-)^+:X^0 \to X^0$ are therefore consistent with the other meaning of the notations $(-)^-$ and $(-)^+$ given in Section~\ref{section:cellular}.

\begin{prop} \label{prop:tauab} (The trace $\tau_{\alpha,\beta}$)
Let \begin{tikzcd}[cramped]f:X\arrow[r,"\sbd"]&Y\end{tikzcd} be a globular subdivision. Let $\alpha,\beta\in Y^0\backslash X^0$. The following statements are equivalent:
\begin{enumerate}
	\item There exists a directed path $\gamma$ of $Y$ from $\alpha$ to $\beta$ such that $\gamma([0,1])\subset |X|\backslash X^0$.
	\item For every cellular decomposition of $X$, $\alpha$ and $\beta$ are in the same globular cell and there exists a unique directed path from $\alpha$ to $\beta$ up to reparametrization (the corresponding trace of $\vec{T}(X)(\alpha,\beta)=\vec{T}(Y)(\alpha,\beta)$ is denoted by $\tau_{\alpha,\beta}$).
\end{enumerate}
\end{prop}

\bpf
By Theorem~\ref{thm:T-homotopy}, and the map $f:X\to Y$ being a globular subdivision, there is the isomorphism of directed spaces $\cont(f):\cont(X)\iso \cont(Y)$, which implies $\vec{T}(X)(\alpha,\beta)=\vec{T}(Y)(\alpha,\beta)$. The implication $(2)\Rightarrow (1)$ is obvious. If $(1)$ holds, then $\alpha$ and $\beta$ belongs by to the same globular cell of $X$; if $\widehat{g}$ is the attaching map, the existence of $\gamma$ implies that $\widehat{g}(z,t)=\alpha$ and $\widehat{g}(z,u)=\beta$ for some $t\leq u$. Hence we obtain $(2)$. 
\epf

The proof of Theorem~\ref{thm:plus-minus-map} relies on \cite[Proposition~8.7 and 8.8 and 8.9 and 8.12]{GlobularNaturalSystem} which are now recalled to facilitate the reading of the proof.

\begin{prop} \label{prop:recall}
Let $X$ be a cellular multipointed $d$-space. Let $c$ be a globular cell of $X$ of dimension greater than $1$. One has
\begin{enumerate}
	\item Let $\alpha\in X^0$ and $x\in c$. There is the homeomorphism $\vec{T}(c^+,\alpha)\iso \vec{T}(x,\alpha)$.
	\item Let $\alpha\in X^0$ and $x\in c$. There is the homeomorphism $\vec{T}(\alpha,c^-)\iso \vec{T}(\alpha,x)$.
	\item Let $\gamma$ be a directed path of $X$ with $\gamma(0),\gamma(1)\in c$. There is the homeomorphism $\{\tau_{\gamma(0),\gamma(1)}\} \sqcup \vec{T}(X)(c^+,c^-) \iso \vec{T}(X)(\gamma(0),\gamma(1))$ ($\tau_{\gamma(0),\gamma(1)}$ is defined in Proposition~\ref{prop:tauab}).
	\item Let $d$ be another globular cell of $X$ of dimension greater than $1$ with $c\neq d$. Let $x\in c$ and $y\in d$. There is the homeomorphism $\vec{T}(X)(c^+,d^-) \iso \vec{T}(X)(x,y)$.
\end{enumerate}
\end{prop}

\bpf $(1)$ is \cite[Proposition~8.7]{GlobularNaturalSystem}. $(2)$ is \cite[Proposition~8.8]{GlobularNaturalSystem}. $(3)$ is \cite[Proposition~8.12]{GlobularNaturalSystem}. Finally, $(4)$ is \cite[Proposition~8.9]{GlobularNaturalSystem}.
\epf

The geometric meaning of $(1)$ is that a trace from $x$ to $\alpha$ has to leave $c$ in a unique way by $c^+$ before reaching $\alpha$ which can be equal to $c^+$. For symmetric reasons, we have $(2)$. The case of $(3)$ is a little bit different. The space of traces from $\gamma(0)$ to $\gamma(1)$ which belong to the same globular cell contains two connected components: the unique trace $\tau_{\gamma(0),\gamma(1)}$ going from $\gamma(0)$ to $\gamma(1)$ without leaving the globular cell $c$ and all other traces leaving $c$ by the unique trace going from $\gamma(0)$ to $c^+$ and returning to $c$ by $c^-$ before reaching $\gamma(1)$ by the unique trace going from $c^-$ to $\gamma(1)$. Finally the last case means geometrically that any trace from $x\in c$ to $y\in d$ with $c\neq d$ has to leave $c$ by the unique trace going from $x$ to $c^+$ and to reach $y$ by reaching at first $d^-$ and by following the unique trace going from $d^-$ to $y$.

\bth \label{thm:plus-minus-map}
Let \begin{tikzcd}[cramped]f:X\arrow[r,"\sbd"]&Y\end{tikzcd} be a globular subdivision. Let $\alpha,\beta\in Y^0$. Then either there exists a directed path $\gamma$ of $Y$ from $\alpha$ to $\beta$ such that $\gamma([0,1])\subset |X|\backslash X^0$ and there is the homeomorphism (see Proposition~\ref{prop:tauab} for the definition of $\tau_{\alpha,\beta}$) \[\vec{T}(Y)(\alpha,\beta) \iso \{\tau_{\alpha,\beta}\} \cup \vec{T}(X)(\alpha^+,\beta^-),\] or there is no such directed path and there is the homeomorphism \[\vec{T}(Y)(\alpha,\beta) \iso \vec{T}(X)(\alpha^+,\beta^-).\]
\eth

\bpf
Choose a cellular decomposition of $X$. Assume that there is a directed path $\gamma$ of $\cont(Y)$ such that $\gamma([0,1])\subset |X|\backslash X^0$. Then $\alpha=\gamma(0)$ and $\beta=\gamma(1)$ belongs to the same globular cell of $X$ of dimension greater than $1$. The homeomorphism $\vec{T}(Y)(\alpha,\beta) \iso \{\tau_{\alpha,\beta}\} \cup \vec{T}(Y)(\alpha^+,\beta^-)$ is a consequence of Proposition~\ref{prop:recall}~$(3)$. Now assume that there is not such a directed path. There are five mutually exclusive possibilities. 
\begin{enumerate}[leftmargin=*]
	\item $\alpha,\beta\in X^0$: $\alpha=\alpha^+$ and $\beta=\beta^-$ implies the equality $\vec{T}(Y)(\alpha,\beta) = \vec{T}(Y)(\alpha^+,\beta^-)$.
	\item $\alpha\in X^0$ and $\beta\in |X|\backslash X^0$: $\alpha=\alpha^+$ and Proposition~\ref{prop:recall}~$(2)$ implies the homeomorphism $\vec{T}(Y)(\alpha,\beta) \iso \vec{T}(Y)(\alpha^+,\beta^-)$.
	\item $\alpha\in |X|\backslash X^0$ and $\beta\in X^0$: $\beta=\beta^-$ and Proposition~\ref{prop:recall}~$(1)$ implies the homeomorphism $\vec{T}(Y)(\alpha,\beta) \iso \vec{T}(Y)(\alpha^+,\beta^-)$.
	\item $\alpha,\beta\in |X|\backslash X^0$ and belonging to two different globular cells of $X$: Proposition~\ref{prop:recall}~$(4)$ implies the homeomorphism $\vec{T}(Y)(\alpha,\beta) \iso \vec{T}(Y)(\alpha^+,\beta^-)$.
	\item $\alpha,\beta\in |X|\backslash X^0$ and belonging to the same globular cell $c$ of $X$: it is a variant of Proposition~\ref{prop:recall}~$(4)$ whose proof is left to the reader; intuitively, every trace from $\alpha$ to $\beta$ has to go out from $c$ by the unique trace $\tau_{\alpha,c^+}$ going from $\alpha$ to $c^+$ and to return to $c$ by $c^-$ followed by the unique trace $\tau_{c^-,\beta}$ from $c^-$ to $\beta$. 
\end{enumerate}
The proof is complete thanks to Proposition~\ref{prop:sbd-path}.
\epf

\begin{figure}
	\[
\def\n{3}
\begin{tikzpicture}
	\draw[-] [thick] (0,0) -- (\n,0);
	\draw[-] [dark-red] [thick] (\n,0) -- (\n,\n);
	\draw[-] [thick] (\n,\n) -- (0,\n);
	\draw[-] [dark-red] [thick] (0,\n) -- (0,0);
	\draw[-] [thick] (\n+2,0) -- (2+2*\n,0);
	\draw[-] [dark-red] [thick] (2+2*\n,0) -- (2+2*\n,\n);
	\draw[-] [thick] (2+2*\n,\n) -- (2+\n,\n);
	\draw[-] [dark-red] [thick] (2+\n,\n) -- (2+\n,0);
	\draw[-] [thick] (0,\n/3) -- (\n,\n/3);
	\draw[-] [dashed] [thick] (0,2*\n/3) -- (\n,2*\n/3);
	\draw[-] [dashed] [thick] (2+\n,\n/3) -- (2+2*\n,\n/3);
	\draw[-] [thick] (2+\n,2*\n/3) -- (2+2*\n,2*\n/3);
\end{tikzpicture}
\]
\caption{Two globular subdivisions of $\globM([0,1])$}
\label{fig:two-gl-subd-same-dspace}
\end{figure}
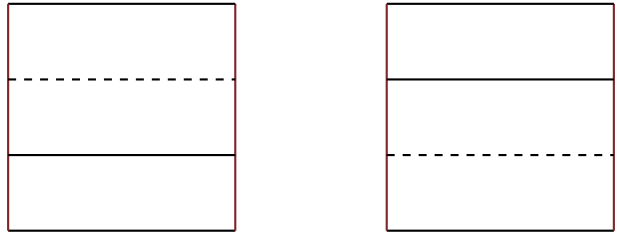

A globular subdivision \begin{tikzcd}[cramped]f:X\arrow[r,"\sbd"]&Y\end{tikzcd} is not necessarily well-behaved with respect to the chosen cellular decompositions of $X$ and $Y$. More precisely, the image $f(\widehat{c})$ of the closure $\widehat{c}$ of a globular cell $c$ of $X$ is not necessarily a cell subcomplex of $Y$. To see that, let us start from the homeomorphism $[0,1/3] \sqcup_{1/3} [1/3,1] \iso [0,2/3] \sqcup_{2/3} [2/3,1]$. This gives rise to two cellular decompositions of $\globM([0,1])$ which can be depicted as in Figure~\ref{fig:two-gl-subd-same-dspace} (it is understood that each red line must be identified to a point). The identity of $\globM([0,1])$ is a globular subdivision which is such a pathological example. 

Theorem~\ref{thm:sbd-cell} proves that it is always possible for any globular subdivision \begin{tikzcd}[cramped]f:X\arrow[r,"\sbd"]&Y\end{tikzcd} to choose another cellular decomposition of $Y$ to avoid this kind of issue.

\bth \label{thm:sbd-cell}
Let \begin{tikzcd}[cramped]f:X\arrow[r,"\sbd"]&Y\end{tikzcd} be a globular subdivision. Let $\widetilde{X}:\lambda\to \ptop{\mathcal{M}}$ be a cellular decomposition of $X$. There exists a transfinite tower of cellular multipointed $d$-spaces $\widetilde{Y}:\lambda\to \ptop{\mathcal{M}}$ and a map of transfinite towers $\widetilde{X}\to \widetilde{Y}$ such that the colimit is the globular subdivision $X\to Y$ and such that for all $\nu<\lambda$, there is the 
commutative diagrams of multipointed $d$-spaces of the form
\[\begin{tikzcd}[row sep=3em, column sep=3em]
	\globM(\mathbf{S}^{n_\nu-1}) \arrow[d,rightarrowtail] \arrow[r,"g_\nu"] \arrow[dd,bend right=60pt,rightarrowtail]& \widetilde{X}_\nu \arrow[d,rightarrowtail] \arrow[r,"\sbd"] & \widetilde{Y}_\nu \arrow[dd,rightarrowtail]\\
	\globM(\mathbf{D}^{n_\nu}) \arrow[d,"\sbd"]\arrow[r,"\widehat{g_\nu}"] & \cocartesian \widetilde{X}_{\nu+1}\arrow[rd,shorten >=0.5em,"\sbd"]  \\
	\globM(\mathbf{D}^{n_\nu})_{F(c_\nu)} \arrow[rr] & & \widetilde{Y}_{\nu+1} \arrow[lu, phantom, "\ulcorner"{font=\Large}, pos=-0.1]
\end{tikzcd}\]
Moreover, the connection maps of the globular subdivision \begin{tikzcd}[cramped]f:X_\nu\arrow[r,"\sbd"]&Y_\nu\end{tikzcd} are the restrictions to $Y_\nu^0$ of the connection maps of the globular subdivision \begin{tikzcd}[cramped]f:X_{\nu+1}\arrow[r,"\sbd"]&Y_{\nu+1}\end{tikzcd} for all $\nu<\lambda$.
\eth

Note that each map $\widehat{Y}_\nu\to \widehat{Y}_{\nu+1}$ is a finite composition of pushouts of maps of $I^{gl,top}\cup\{C\}$ by Proposition~\ref{prop:raf-cell}.

The hypothesis that $Y$ is cellular is used to guarantee that $Y^0$ is discrete and at the very end of the proof for using Theorem~\ref{thm:T-homotopy}. Without this hypothesis, i.e. by only assuming that $f$ induces a homeomorphism, and even by assuming $Y^0$ discrete, there is no guarantee for the map $\widetilde{Y}_\lambda\to Y$, which induces a homeomorphism between the set of states and between the underlying spaces, to be an isomorphism: the multipointed $d$-space $Y$  could have \textit{more} execution paths than $\widetilde{Y}_\lambda$ indeed.

\bpf 
There is a bijection of sets 
\[
|X| = \coprod_{c\in \mathcal{C}(X)} c
\]
where $\mathcal{C}(X)$ is the set of globular cells of $\widetilde{X}$. We obtain a bijection of sets 
\[
|Y| = \coprod_{c\in \mathcal{C}(X)} f(c).
\]
Each cell $c_\nu$ of $X$ for $\nu<\lambda$ corresponds to a pushout diagram of the form 
\[\begin{tikzcd}[row sep=3em, column sep=3em]
	\globM(\mathbf{S}^{n_\nu-1}) \arrow[d,rightarrowtail] \arrow[r,"g_\nu"] & \widetilde{X}_\nu \arrow[d,rightarrowtail] \\
	\globM(\mathbf{D}^{n_\nu}) \arrow[r,"\widehat{g_\nu}"] & \cocartesian \widetilde{X}_{\nu+1}
\end{tikzcd}\]
for some $n_\nu\geq 0$. Consider a globular cell $c$ of $X$ with $\dim(c)\geq 1$. The closure $\widehat{f(c)}$ of each $f(c)$ is a compact in the Hausdorff space $|X|=|Y|$. Therefore, the set $f(c) \cap Y^0 \subset \widehat{f(c)}\cap Y^0$ is finite, $Y^0$ being discrete because $Y$ is cellular by hypothesis. Let $F(c)= f(c) \cap Y^0$ for $c$ running over the set $\mathcal{C}(X)$ of globular cells of $X$ of dimension greater than $1$. We are going to construct by transfinite induction a transfinite tower $\widetilde{Y}:\lambda\to \ptop{\mathcal{M}}$ and a map of transfinite towers $\widetilde{X}\to \widetilde{Y}$ such that for all $\nu\leq \lambda$, there is a homeomorphism $\widetilde{X}_\nu\iso \widetilde{Y}_\nu$ as follows. Let $\widetilde{Y}_0=X^0$. Assume that the map of transfinite towers $\widetilde{X}\to \widetilde{Y}$ is constructed until some ordinal $\nu<\lambda$. Consider the commutative diagram of solid arrows of multipointed $d$-spaces 
\[\begin{tikzcd}[row sep=3em, column sep=3em]
	\globM(\mathbf{S}^{n_\nu-1}) \arrow[dd,rightarrowtail]\arrow[dr,equal] \arrow[rr,"g_\nu"]&& \widetilde{X}_\nu \arrow[dd,rightarrowtail] \arrow[rd]\\
	& \globM(\mathbf{S}^{n_\nu-1}) \arrow[rr,crossing over]  && \widetilde{Y}_\nu \arrow[dd,rightarrowtail]\\
	\globM(\mathbf{D}^{n_\nu}) \arrow[rr,"\widehat{g_{\nu}}",pos=0.3]\arrow[rd,"\sbd"] && \cocartesian \widetilde{X}_{\nu+1} \arrow[rd,dashed,shorten >=0.5em]\\
	& \globM(\mathbf{D}^{n_\nu})_{F(c_\nu)} \arrow[uu,leftarrowtail,crossing over]\arrow[rr] && \arrow[lu, phantom, "\ulcorner"{font=\Large}, pos=-0.1] \widetilde{Y}_{\nu+1}
\end{tikzcd}\]
%depicted in Figure~\ref{fig:passage}. 
The universal property of the pushout yields a (unique) map $\widetilde{X}_{\nu+1}\to \widetilde{Y}_{\nu+1}$. Since the functor $X\mapsto (|X|,X^0)$ from $\ptop{\mathcal{M}}$ to $\mtop$ is topological, it is colimit-preserving. Therefore, we obtain from the homeomorphism $|\widetilde{X}_\nu| \iso |\widetilde{Y}_\nu|$ the homeomorphism $|\widetilde{X}_{\nu+1}| \iso |\widetilde{Y}_{\nu+1}|$. For a limit ordinal $\nu\leq \lambda$, let us define the map of multipointed $d$-spaces $\widetilde{X}_{\nu}\to \widetilde{Y}_{\nu}$ as a colimit. Since the functor $X\mapsto |X|$ is colimit-preserving, we obtain the homeomorphism $|\widetilde{X}_\nu| \iso |\widetilde{Y}_\nu|$ also for the limit ordinals $\nu\leq \lambda$. By Proposition~\ref{prop:raf-cell}, each multipointed $d$-space $\widetilde{Y}_\nu$ for $\nu\leq \lambda$ is cellular. We deduce that the maps $\widetilde{X}_\nu\to \widetilde{Y}_\nu$ are globular subdivisions for all $\nu\leq \lambda$. By reorganizing the cube above, we obtain the planar commutative diagram of multipointed $d$-spaces of the statement of the theorem for all $\nu<\lambda$. The map of multipointed $d$-spaces $\widetilde{X}_0 \to Y$ factors uniquely as a composite $\widetilde{X}_0\iso \widetilde{Y}_0 \to Y$. Assume that for an ordinal $\nu<\lambda$, the map of multipointed $d$-spaces $\widetilde{X}_\nu \to Y$ factors uniquely as a composite $\widetilde{X}_\nu\iso \widetilde{Y}_\nu \to Y$. From the commutative square of multipointed $d$-spaces 
\[
\begin{tikzcd}[row sep=4em,column sep=4em]
	\globM(\mathbf{S}^{n_\nu-1}) \arrow[r,"{g_\nu}"] \arrow[d]& \widetilde{X}_{\nu}  \arrow[r]& Y\arrow[d,equal]\\
	\globM(\mathbf{D}^{n_\nu})_{F(c_\nu)} \arrow[rr]& &Y
\end{tikzcd}
\]
we deduce that the map of multipointed $d$-spaces $\widetilde{X}_{\nu+1} \to Y$ factors uniquely as a composite $\widetilde{X}_{\nu+1}\iso \widetilde{Y}_{\nu+1} \to Y$. We obtain by transfinite induction that the globular subdivision $X\to Y$ factors uniquely as a composite of globular subdivisions $X\to \widetilde{Y}_\lambda \to Y$. It remains to prove that every execution path of $Y$ is an execution path of $\widetilde{Y}_\lambda$ to complete the proof. Consider an execution path $\gamma$ of $Y$. It is a directed path of the directed space $\cont(Y)$ by definition of the functor $\cont:\ptop{\mathcal{M}} \to \ptop{}$. By Theorem~\ref{thm:T-homotopy}, one has $\cont(X)=\cont(Y)$. Thus, $\gamma$ is a directed path of the directed space $\cont(X)$. By Theorem~\ref{thm:continuous}, $X$ being cellular, there exists an execution path $\gamma'$ of $X$ and $\phi\in \mathcal{I}$ such that $\gamma=\gamma'\phi$. In plain English, every execution path of $Y$ is a piece of an execution path of $X$ between two points of $Y^0$. Therefore it is an execution path of $\widetilde{Y}_\lambda$, which means that $\widetilde{Y}_\lambda = Y$. Finally, the location of the q-cofibrations is a consequence of Proposition~\ref{prop:raf-cell}. \epf

\section{Generating subdivision}
\label{section:glT}

\bd A poset $(P,\leq)$ is \textit{bounded} if there exist $\widehat{0}\in P$ and $\widehat{1}\in P$ such that $P = [\widehat{0},\widehat{1}]$ and such that $\widehat{0} \neq \widehat{1}$. Let $\widehat{0}=\min P$ (the bottom element) and $\widehat{1}=\max P$ (the top element). 
\ed

\begin{nota} \label{definitiondeT} \cite[Definition~4.4]{3eme} \label{notation:Tdef}
	Let $\T$ be the class of inclusions of finite bounded posets $P_1\subset P_2$ preserving the bottom element and the top element. The class $\T$ is essentially small.
\end{nota}

\bd \label{def:generating-sbd}
A \textit{generating subdivision} is a q-cofibration of flows $f^{cof}: P_1^{cof} \to P_2^{cof}$ between q-cofibrant flows such that there exists a commutative square of flows 
\[
\begin{tikzcd}[row sep=3em, column sep=3em]
	P_1^{cof}\arrow[d,rightarrowtail,"f^{cof}"']\arrow[r,"\simeq"]  & P_1 \arrow[d,"f"] \\
	P_2^{cof} \arrow[r,"\simeq"] & \ P_2
\end{tikzcd}
\]
with $f:P_1\subset P_2\in\T$ and such that the horizontal maps are weak equivalences of the q-model structure of flows. 
\ed

In \cite{3eme,4eme,hocont}, such a map is called a generating T-homotopy equivalence. For the same reason as in \cite[Definition~9.1]{GlobularNaturalSystem} where the T-homotopy equivalence terminology is abandoned to the more appropriate globular subdivision terminology, we want to forget the old terminology which was a bit naive.

\begin{nota}
	Let $\T^{cof}$ be an \textit{arbitrary} choice of generating subdivisions $f^{cof}$ for $f$ running over the class of maps $\T$. The class $\T^{cof}$ is essentially small.
\end{nota}

We want to introduce in this section a very specific family $\T^{gl}$ of generating subdivisions which makes the calculations easy for Proposition~\ref{prop:pre-Thomfund1} and Theorem~\ref{thm:plus-minus-map-flow}. We will prove in Proposition~\ref{prop:factorization-gl-cof} that this arbitrary choice has no consequence for the sequel.

\begin{prop} \label{prop:defW}
	Let $n\geq 0$. Consider a finite set $F\subset |\globM(\mathbf{D}^n)| \backslash |\globM(\mathbf{S}^{n-1})|$. There exists a multipointed $d$-space $W^n_F$ and a factorization 
	\[
	\begin{tikzcd}[row sep=3em,column sep=3em]
		\globM(\mathbf{D}^n) \arrow[r,rightarrowtail,"i^n_F"] & W^n_F \arrow[r,twoheadrightarrow,"\simeq"]& \globM(\mathbf{D}^n)_F
	\end{tikzcd}
	\]
	of the globular subdivision $\globM(\mathbf{D}^n)\to \globM(\mathbf{D}^n)_F$ such that 
	\begin{itemize}
		\item One has $i^n_F\in \cell(I^{gl,top}\cup\{C\})$.
		\item There is the isomorphism $W^n_F\rest_{\{0,1\}} \iso \globM(\mathbf{D}^{n+1})$,
		\item The map $\globM(\mathbf{D}^n)\to W^n_F\rest_{\{0,1\}} \iso \globM(\mathbf{D}^{n+1})$ is induced by the inclusion of $\mathbf{D}^n$ into one copy of $\mathbf{D}^n$ in $\mathbf{D}^n\sqcup_{\mathbf{S}^{n-1}}\mathbf{D}^n \iso \mathbf{S}^n\subset \mathbf{D}^{n+1}$. This implies that the map $\globM(\mathbf{D}^n)\to W^n_F\rest_{\{0,1\}}$ is a trivial q-cofibration of multipointed $d$-spaces.
		\item There are the homeomorphisms $\P_{\alpha,\beta}^{top}W^n_F\iso \P_{\alpha,\beta}^{top}\globM(\mathbf{D}^n)_F$ for all $(\alpha,\beta)\neq(0,1)$. 
	\end{itemize}
	Moreover, the map $W^n_F\to \globM(\mathbf{D}^n)_F$ is a weak equivalence between cellular multipointed $d$-spaces. Finally, the map $\dcat(W^n_F)\to \glob(\mathbf{D}^n)_F$ is a trivial q-cofibration between cellular flows.
\end{prop}

Note that the underlying space of $W^n_F$ is $|\globM(\mathbf{D}^{n+1})|$.

\bpf
The multipointed $d$-space $W^n_F$ is obtained in two steps. The first step is the pushout diagram of multipointed $d$-spaces 
\[
\begin{tikzcd}[column sep=3em,row sep=3em]
	\globM(\mathbf{S}^{n-1}) \arrow[d,rightarrowtail,"\subset"'] \arrow[r,"\subset"] & \globM(\mathbf{D}^n)\arrow[d,rightarrowtail]\\
	\globM(\mathbf{D}^n)_F \arrow[r] & \cocartesian V^n_F
\end{tikzcd}
\]
The multipointed $d$-space $V_F^n$ is cellular by Proposition~\ref{prop:raf-cell}. It satisfies \[V_F^n\rest_{\{0,1\}} \iso \globM(\mathbf{D}^{n}\sqcup_{\mathbf{S}^{n-1}} \mathbf{D}^n)\] because $\globM(\mathbf{D}^n)_F\rest_{\{0,1\}}\iso \globM(\mathbf{D}^n)$. Moreover, one has
\[
\forall (\alpha,\beta)\neq(0,1), \P_{\alpha,\beta}^{top}V^n_F \iso \P_{\alpha,\beta}^{top}\globM(\mathbf{D}^n)_F.
\]
After choosing a homeomorphism $\mathbf{D}^{n}\sqcup_{\mathbf{S}^{n-1}} \mathbf{D}^n\iso \mathbf{S}^n$, the second step is the pushout diagram of multipointed $d$-spaces 
\[
\begin{tikzcd}[column sep=3em,row sep=3em]
	\globM(\mathbf{S}^{n}) \iso V_F^n\rest_{\{0,1\}}\arrow[d,rightarrowtail] \arrow[r,"\subset"] & V^n_F\arrow[d,rightarrowtail]\\
	\globM(\mathbf{D}^{n+1}) \arrow[r] & \cocartesian W^n_F
\end{tikzcd}
\]
By Theorem~\ref{thm:colim}, the map $\dcat(i^n_F):\glob(\mathbf{D}^n)\to \dcat(W^n_F)$ belongs to $\cell(I^{gl}\cup\{C\})$, i.e. it is a q-cofibration of flows, the multipointed $d$-space $V^n_F$ being cellular. Since both $W^n_F$ and $\globM(\mathbf{D}^n)_F$ are q-cofibrant, using Theorem~\ref{thm:rappel-all}~$(2)$, the image by the functor $\dcat:\ptop{\mathcal{M}}\to \dtop$ of the weak equivalence $W^n_F\to \globM(\mathbf{D}^n)_F$ is a weak equivalence of the q-model structure of flows between two q-cofibrant flows, and also a q-cofibration between cellular flows by Theorem~\ref{thm:colim}. 
\epf

\begin{prop} \label{prop:Tgl}
	With the notations of Proposition~\ref{prop:poset} and Proposition~\ref{prop:defW}. The q-cofibration \[\dcat(i^n_F):\glob(\mathbf{D}^n)\longrightarrow \dcat(W^n_F)\] is a generating subdivision; more precisely, there exists a commutative diagram of flows 
	\[
	\begin{tikzcd}[row sep=3em, column sep=3em]
		\glob(\mathbf{D}^n)\arrow[d,rightarrowtail,"\dcat(i^n_F)"']\arrow[r,"\simeq",twoheadrightarrow]  & \{0<1\} \arrow[d,"\subset"] \\
		\dcat(W^n_F) \arrow[r,"\simeq",twoheadrightarrow] & (\{0,1\}\cup F,\leq)
	\end{tikzcd}
	\]
	such that the horizontal maps are weak equivalences of the q-model structure of flows and where the map of posets  $\{0<1\}\subset (\{0,1\}\cup F,\leq)$ is defined in Proposition~\ref{prop:poset}.
\end{prop}

\bpf
It is a consequence of Proposition~\ref{prop:poset} and Proposition~\ref{prop:defW}. The two horizontal maps are q-fibrations since the space of execution paths of the flows $\{0<1\}$ and $(\{0,1\}\cup F,\leq)$ are discrete.
\epf

\begin{nota} \label{notation:Igl}
	Consider the set of generating subdivisions \[\T^{gl} = \{\dcat(i^n_F):\glob(\mathbf{D}^n)\to \dcat(W^n_F)\}\] with $n\geq 0$ and $F$ running over all finite subsets of $|\globM(\mathbf{D}^n)| \backslash |\globM(\mathbf{S}^{n-1})|$.
\end{nota}

\begin{prop} \label{prop:pre-Thomfund1}
	Consider a pushout diagram of multipointed $d$-spaces of the form 
	\[
	\begin{tikzcd}[row sep=3em, column sep=3em]
		\globM(\mathbf{D}^n)\arrow[d,rightarrowtail,"i^n_F"']\arrow["g",r]  & A \arrow[d,"f",rightarrowtail] \\
		W^n_F \arrow[r,"\widehat{g}"] & \cocartesian X
	\end{tikzcd}
	\]
	where $A$ is a cellular multipointed $d$-space. Then the map $f$ induces a homotopy equivalence 
	\[
	\II(\dcat(A))(\alpha,\beta) \simeq \II(\dcat(X))(f(\alpha),f(\beta))
	\]
	between q-cofibrant spaces for all $(\alpha,\beta) \in X^0\p X^0$.
\end{prop}

\bpf Consider the pushout diagram of multipointed $d$-spaces 
\[
\begin{tikzcd}[column sep=3em,row sep=3em]
	\globM(\mathbf{D}^n)\arrow[d,rightarrowtail] \arrow["g",r] & A \arrow[d,rightarrowtail]\arrow[dd,bend left=60pt,"f",rightarrowtail]\\
	\globM(\mathbf{D}^{n+1})\iso W^n_F\rest_{\{0,1\}} \arrow[r]\arrow[d,"\sbd"']& Z \arrow[d,"\sbd"'] \cocartesian \\
	W^n_F \arrow[r,"\widehat{g}"] & \cocartesian X
\end{tikzcd}
\]
The map $Z\to X$ is a globular subdivision, $Z$ being cellular and the underlying space functor being colimit-preserving by Proposition~\ref{prop:underlyingspace-almost-left-Quillen-adjoint}. We deduce the isomorphism of multipointed $d$-spaces $Z\iso X\rest_{A^0}$. We obtain the commutative diagram of multipointed $d$-spaces
\[
\begin{tikzcd}[column sep=3em,row sep=3em]
	\globM(\mathbf{D}^n)\arrow[d,rightarrowtail] \arrow["g",r] & A \arrow[d,rightarrowtail]\arrow[dd,bend left=60pt,"f",rightarrowtail]\\
	\globM(\mathbf{D}^{n+1})\iso W^n_F\rest_{\{0,1\}} \arrow[r]\arrow[d,"\sbd"']& X\rest_{A^0} \arrow[d,"\sbd"'] \cocartesian \\
	W^n_F \arrow[r,"\widehat{g}"] & \cocartesian X
\end{tikzcd}
\]
By Theorem~\ref{thm:colim} , the image by $\dcat$ of the top commutative square yields the pushout diagram of flows 
\[
\begin{tikzcd}[column sep=3em,row sep=3em]
	\glob(\mathbf{D}^n)\arrow[d,rightarrowtail,"\simeq"] \arrow["\dcat(g)",r] & \dcat(A) \arrow[d,rightarrowtail,"\simeq"]\\
	\glob(\mathbf{D}^{n+1}) \arrow[r]& \dcat(X\rest_{A^0})  \cocartesian
\end{tikzcd}
\]
such that the vertical maps are q-cofibrations of flows, which moreover are trivial since the inclusion $\mathbf{D}^{n}\subset \mathbf{D}^{n+1}$ is a trivial q-cofibration of spaces. By Theorem~\ref{thm:complement-flow}, for all $(\alpha,\beta) \in X^0\p X^0$, the map $\P_{\alpha,\beta}\dcat(A)\to \P_{\alpha,\beta}\dcat(X\rest_{A^0})$ is then a trivial q-cofibration of spaces between q-cofibrant spaces. Hence it is also a homotopy equivalence. We therefore obtain for all $(\alpha,\beta) \in X^0\p X^0$ that 
\[
\II(\dcat(A))(\alpha,\beta) \simeq \II(\dcat(X\rest_{A^0}))(\alpha,\beta) \iso \II(\dcat(X))(f(\alpha),f(\beta)),
\]
the homotopy equivalence by definition of the functor $\II$ (see Notation~\ref{def:II}), and the homeomorphism by Proposition~\ref{prop:sbd-path}, the map $X\rest_{A^0}\to X$ being a globular subdivision. 
\epf

Theorem~\ref{thm:plus-minus-map-flow} is the analogue for the maps of $\T^{gl}$ of Theorem~\ref{thm:plus-minus-map} for the globular subdivisions.

\bth \label{thm:plus-minus-map-flow}
Consider a pushout diagram of multipointed $d$-spaces of the form 
\[
\begin{tikzcd}[row sep=3em, column sep=3em]
	\globM(\mathbf{D}^n)\arrow[d,rightarrowtail,"i^n_F"']\arrow["g",r]  & A \arrow[d,"f",rightarrowtail] \\
	W^n_F \arrow[r,"\widehat{g}"] & \cocartesian X
\end{tikzcd}
\]
where $A$ is a cellular multipointed $d$-space. Consider the set maps $(-)^+,(-)^-:X^0 \to A^0$ defined by $\alpha^+=\alpha^-=\alpha$ if $\alpha\in A^0$ and $\alpha^+=g(1)=\widehat{g}(1)$ and $\alpha^-=g(0)=\widehat{g}(0)$ if $\alpha\in X^0\backslash A^0$. Let $\alpha,\beta\in X^0$. There are two mutually exclusive cases (see Proposition~\ref{prop:tauab} for the definition of $\tau_{\alpha,\beta}$):
\begin{enumerate}
	\item $\alpha,\beta$ in $X^0\backslash A^0$ and $\II(\dcat(W^n_F))(\alpha,\beta)\neq \varnothing$; in this case, there is the homotopy equivalence \[\II(\dcat(X))(\alpha,\beta) \simeq \{\tau_{\alpha,\beta}\} \sqcup \II(\dcat(A))(\alpha^+,\beta^-).\]
	\item There is the homotopy equivalence \[\II(\dcat(X))(\alpha,\beta) \simeq \II(\dcat(A))(\alpha^+,\beta^-)\] otherwise.
\end{enumerate}
\eth

\bpf
By Proposition~\ref{prop:raf-cell} and Theorem~\ref{thm:colim}, there is the pushout diagram of cellular flows
\[
\begin{tikzcd}[row sep=3em, column sep=3em]
	\glob(\mathbf{D}^n)\arrow[d,rightarrowtail,"\dcat(i^n_F)"']\arrow["\dcat(g)",r]  & \dcat(A) \arrow[d,"f",rightarrowtail] \\
	\dcat(W^n_F) \arrow[r,"\dcat(\widehat{g})"] & \cocartesian \dcat(X)
\end{tikzcd}
\]
The functor $\II:\dtop\to \cat_\top$ being a left adjoint, there is the pushout diagram of enriched small categories 
\[
\begin{tikzcd}[row sep=3em, column sep=3em]
	\II(\glob(\mathbf{D}^n))\arrow[d,"\II(i^n_F)"']\arrow["\II(g)",r]  & \II(\dcat(A)) \arrow[d,"\II(f)"] \\
	\II(\dcat(W^n_F)) \arrow[r,"\II(\widehat{g})"] & \cocartesian \II(\dcat(X))
\end{tikzcd}
\]
Assume first that $\alpha,\beta\in X^0\backslash A^0$ and that $\P_{\alpha,\beta}^{top}W^n_F\neq \varnothing$ or $\alpha=\beta$, i.e. $\II(\dcat(W^n_F))(\alpha,\beta)$ is not empty. By Proposition~\ref{prop:defW} (i.e. by construction of the multipointed $d$-space $W^n_F$), there is the homeomorphism
\begin{multline*}
	\II(\dcat(X))(\alpha,\beta)\iso \II(\dcat(W^n_F))(\alpha,\beta) \sqcup\\ \II(\dcat(W^n_F))(\alpha,g(1))\p \II(\dcat(X))(g(1),g(0)) \p \II(\dcat(W^n_F))(g(0),\beta).
\end{multline*} 
Since the map $\dcat(W^n_F)\to \glob(\mathbf{D}^n)_F$ is a weak equivalence of the q-model structure of flows between q-cofibrant flows by Proposition~\ref{prop:defW}, there are the homotopy equivalences 
\begin{align*}
	& \II(\dcat(W^n_F))(\alpha,\beta) \simeq \II(\glob(\mathbf{D}^n_F))(\alpha,\beta)\\
	& \II(\dcat(W^n_F))(\alpha,g(1)) \simeq \II(\glob(\mathbf{D}^n_F)))(\alpha,g(1))\\
	& \II(\dcat(W^n_F))(g(0),\beta) \simeq \II(\glob(\mathbf{D}^n_F))(g(0),\beta)
\end{align*}
by Theorem~\ref{thm:carac-weakeq-trace}. Thus the topological spaces $\II(\dcat(W^n_F))(\alpha,\beta)$, $\II(\dcat(W^n_F))(\alpha,g(1))$ and $\II(\dcat(W^n_F))(g(0),\beta)$ are contractible. Using Proposition~\ref{prop:pre-Thomfund1}, we obtain the homotopy equivalence 
\[
\II(\dcat(X))(\alpha,\beta) \simeq \{\tau_{\alpha,\beta}\} \sqcup \II(\dcat(A))(g(1),g(0)) \iso \{\tau_{\alpha,\beta}\} \sqcup \II(\dcat(A))(\alpha^+,\beta^-).
\]
The case $\alpha,\beta\in X^0\backslash A^0$ and $\II(\dcat(W^n_F))(\alpha,\beta) =\varnothing$ leads to the homotopy equivalence 
\[
\II(\dcat(X))(\alpha,\beta) \simeq  \II(\dcat(A))(g(1),g(0)) \iso \II(\dcat(A))(\alpha^+,\beta^-) 
\]
for similar reasons. The other cases are similar to the one treated in the proof of Theorem~\ref{thm:plus-minus-map}. It remains three mutually exclusive cases. 
\begin{enumerate}[leftmargin=*]
	\item $\alpha,\beta\in A^0$: Proposition~\ref{prop:pre-Thomfund1} implies the homotopy equivalence \[\II(\dcat(X))(\alpha,\beta) \simeq \II(\dcat(A))(\alpha,\beta).\]
	\item $\alpha\in A^0$ and $\beta\in X^0\backslash A^0$: by Proposition~\ref{prop:defW}, we obtain the homeomorphism \[\II(\dcat(X))(\alpha,\beta) \iso \II(\dcat(X))(\alpha,g(0)))\p \II(\dcat(W^n_F))(g(0),\beta),\] and by Proposition~\ref{prop:pre-Thomfund1} the homotopy equivalence \[\II(\dcat(X))(\alpha,\beta) \simeq \II(\dcat(A))(\alpha,g(0)),\] the space $\II(\dcat(W^n_F))(g(0),\beta)$ being contractible by Proposition~\ref{prop:defW}.
	\item $\alpha\in X^0\backslash A^0$ and $\beta\in A^0$: by Proposition~\ref{prop:defW}, we obtain the homeomorphism \[\II(\dcat(X))(\alpha,\beta) \iso \II(\dcat(W^n_F))(\alpha,g(1)))\p \II(\dcat(X))(g(1),\beta),\] and by Proposition~\ref{prop:pre-Thomfund1} the homotopy equivalence \[\II(\dcat(X))(\alpha,\beta) \simeq \II(\dcat(A))(g(1),\beta),\] the space $\II(\dcat(W^n_F))(\alpha,g(1)))$ being contractible by Proposition~\ref{prop:defW}.
\end{enumerate}
\epf

\section{Globular subdivision and generating subdivision}
\label{section:final}

Proposition~\ref{prop:factorization-gl-cof} is required before proving Theorem~\ref{thm:transfinite-sequence-tower}. 

\begin{prop} \label{prop:factorization-gl-cof}
	Let $n\geq 0$. Let $F$ be a finite subset of $|\globM(\mathbf{D}^n)|\backslash |\globM(\mathbf{S}^{n-1})|$. Consider a pushout diagram of flows 
	\[
	\begin{tikzcd}[row sep=3em,column sep=3em]
		\glob(\mathbf{D}^n) \arrow[d,rightarrowtail,"\dcat(i^n_F)"']\arrow[r] & A \arrow[d,"f",rightarrowtail]\\
		\dcat(W^n_F) \arrow[r] & \cocartesian X
	\end{tikzcd}
	\]
	with $A$ and $X$ q-cofibrant. Then the q-cofibration $f:A\to X$ factors as a composite 
	\[
	\begin{tikzcd}[column sep=5em,row sep=3em]
		f:A \arrow[r,rightarrowtail,"\in \cell(\T^{cof})"] & X' \arrow[r,"\simeq"]& X
	\end{tikzcd}
	\] 
	of a map of $\cell(\T^{cof})$ followed by a weak equivalence of flows.
\end{prop}

\bpf
By Proposition~\ref{prop:Tgl}, there is a commutative diagram of solid arrows of q-cofibrant flows
\[
\begin{tikzcd}[row sep=3em, column sep=5em]
	&&\\
	\{0<1\}^{cof} \arrow[d,rightarrowtail,"\in \T^{cof}"] \arrow[dashed,r,"h_1"] \arrow[rr,bend left=20pt,"\simeq"]  & \glob(\mathbf{D}^n) \arrow[d,rightarrowtail,"\dcat(i^n_F)"] \arrow[r,twoheadrightarrow,"\simeq"] & \{0<1\} \arrow[d] \\
	(\{0,1\}\cup F,\leq)^{cof} \arrow[dashed,r,"h_2"] \arrow[rr,bend right=20pt,"\simeq"'] & \dcat(W^n_F) \arrow[twoheadrightarrow,r,"\simeq"] & (\{0,1\}\cup F,\leq)\\
	&&
\end{tikzcd}
\]
Using Proposition~\ref{prop:useful-lemma}, we deduce the existence of $h_1$ and $h_2$ making the diagram above commutative. By the \ttt, the maps $h_1$ and $h_2$ are weak equivalences of flows. The commutative diagram of flows above can then be reorganized as the cube 
\[\begin{tikzcd}[row sep=3em, column sep=3em]
	\{0<1\}^{cof} \arrow[dd,rightarrowtail]\arrow[dr,"h_1"] \arrow[rr]&& A \arrow[dd,rightarrowtail] \arrow[rd,equal]\\
	& \glob(\mathbf{D}^n) \arrow[rr,crossing over]  && A \arrow[dd,rightarrowtail]\\
	(\{0,1\}\cup F,\leq)^{cof} \arrow[rr]\arrow[rd,"h_2"] && \cocartesian X' \arrow[rd,dashed,"s",shorten >=0.5em]\\
	& \dcat(W^n_F) \arrow[uu,"\dcat(i^n_F)",leftarrowtail,crossing over,pos=0.8]\arrow[rr] && \arrow[lu, phantom, "\ulcorner"{font=\Large}, pos=-0.1] X
\end{tikzcd}\] 
where the existence of the flow $X'$ and  of the map $s:X'\to X$ comes from the universal property of the pushout. By the cube lemma \cite[Lemma~5.2.6]{MR99h:55031}, the map $s$ is a weak equivalence of flows.
\epf

\bth \label{thm:transfinite-sequence-tower}
Let $n\geq 0$. Let $F$ be a finite subset of $|\globM(\mathbf{D}^{n})|\backslash |\globM(\mathbf{S}^{n-1})|$. Consider a commutative diagram of cellular multipointed $d$-spaces 
\[\begin{tikzcd}[row sep=3em, column sep=3em]
	\globM(\mathbf{S}^{n-1}) \arrow[dd,bend right=60pt,rightarrowtail]\arrow[d,rightarrowtail] \arrow[r,"g"] & X \arrow[d,rightarrowtail] \arrow[r,"f"] & Y \arrow[dd,rightarrowtail]\\
	\globM(\mathbf{D}^{n}) \arrow[d,"\sbd"]\arrow[r,"\widehat{g}"] & \cocartesian \overline{X}\arrow[rd,shorten >=0.5em,"\overline{f}"]  \\
	\globM(\mathbf{D}^{n})_{F} \arrow[rr] & & \overline{Y} \arrow[lu, phantom, "\ulcorner"{font=\Large}, pos=-0.1]
\end{tikzcd}\]
such that $f:X\to Y$ and $\overline{f}:\overline{X}\to \overline{Y}$ are globular subdivisions. Assume that $\dcat(f)$ factors as a composite $\dcat(f)=w.i$ where $i:\dcat(X)\to Z$ belongs to $\cell(\T^{cof})$ and where $w:Z\to\dcat(Y)$ is a weak equivalence of the q-model structure of flows for some flow $Z$. Then there exists a commutative diagram of flows
\[\begin{tikzcd}[row sep=3em, column sep=3em]
	\dcat(X) \arrow[rr,bend left=30pt,"\dcat(f)"] \arrow[d,rightarrowtail] \arrow[r,"i",rightarrowtail] & Z \arrow[r,"w"] \arrow[d,rightarrowtail]& \dcat(Y) \arrow[d,rightarrowtail]\\
	\dcat(\overline{X}) \arrow[rr,bend right=30pt,"\dcat(\overline{f})"'] \arrow[r,"\overline{\imath}",rightarrowtail] & \overline{Z} \arrow[r,"\overline{w}"] & \dcat(\overline{Y})
\end{tikzcd}\]
such that $\overline{\imath}:\dcat(\overline{X})\to \overline{Z}$ belongs to $\cell(\T^{cof})$ and  $\overline{w}:\overline{Z}\to \dcat(\overline{Y})$ is a weak equivalence of the q-model structure of flows, and finally such that the canonical map \[U=\dcat(\overline{X})\sqcup_{\dcat(X)} Z \to \overline{Z}\] induced by the universal property of the pushout belongs to $\cell(\T^{cof})$.
\eth

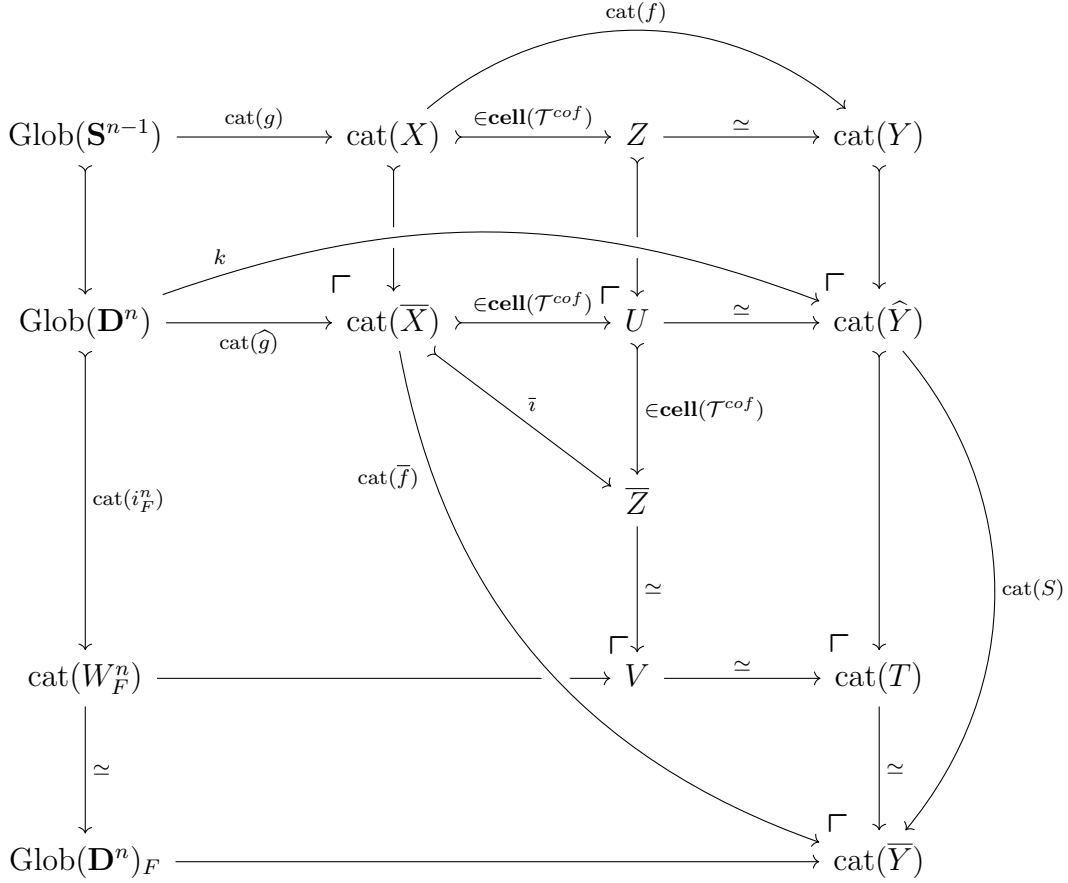
\begin{figure}
	\[
	\begin{tikzcd}[row sep=4em, column sep=5em]
		\glob(\mathbf{S}^{n-1})\arrow[r,"\dcat(g)"]\arrow[rightarrowtail,d]& \dcat(X) \arrow[rr,bend left=40pt,"\dcat(f)"]\arrow[r,rightarrowtail,"\in\cell(\T^{cof})"] \arrow[rightarrowtail,d] & Z  \arrow[r,"\simeq"]\arrow[rightarrowtail,d] & \arrow[rightarrowtail,d] \dcat(Y) \\
		\glob(\mathbf{D}^{n})\arrow[rrr,bend left=20pt,crossing over,"k",pos=0.1]\arrow[dd,rightarrowtail,"\dcat(i^n_F)"]\arrow[r,"\dcat(\widehat{g})"'] & \cocartesian \dcat(\overline{X}) \arrow[rd,"\overline{\imath}",rightarrowtail] \arrow[r,rightarrowtail,"\in\cell(\T^{cof})"] & \cocartesian U  \arrow[r,"\simeq"]\arrow[d,rightarrowtail,"\in\cell(\T^{cof})"] & \cocartesian \dcat(\widehat{Y}) \arrow[dd,rightarrowtail] \arrow[ddd,bend left=40pt,"\dcat(S)"]\\
		&& \overline{Z}\arrow[d,"\simeq"]\\
		\dcat(W^n_F) \arrow[d,"\simeq"]\arrow[rr]&& \arrow[luu, phantom, "\ulcorner"{font=\Large}, pos=0] V \arrow[r,"\simeq"] & \cocartesian \dcat(T) \arrow[d,"\simeq"]\\
		\glob(\mathbf{D}^{n})_F \arrow[rrr]&&& \arrow[lluuu,leftarrow,bend left=30pt,"{\dcat(\overline{f})}",pos=0.85,crossing over]\dcat(\overline{Y}) \cocartesian
	\end{tikzcd}
	\]
	\caption{Adding and subdividing a globe $\glob(\mathbf{D}^{n})$ to $\dcat(X)$}
	\label{fig:add-and-subdivide-one-globe}
\end{figure}

\bpf
Let $\widehat{Y} = \overline{X}\sqcup_X Y$. We obtain the commutative diagram of cellular multipointed $d$-spaces 
\[\begin{tikzcd}[row sep=3em, column sep=3em]
	\globM(\mathbf{S}^{n-1}) \arrow[d] \arrow[r,"g"] & X \arrow[d,rightarrowtail] \arrow[r,"\sbd"] & Y \arrow[d,rightarrowtail]\\
	\globM(\mathbf{D}^{n}) \arrow[d,"\sbd"]\arrow[r,"\widehat{g}"] & \cocartesian \overline{X}\arrow[r,"\sbd"]\arrow[rd,shorten >=0.5em,"\sbd"] & \cocartesian\widehat{Y} \arrow[d,"S"]\\
	\globM(\mathbf{D}^{n})_{F} \arrow[rr] & & \overline{Y} \arrow[lu, phantom, "\ulcorner"{font=\Large}, pos=-0.05]
\end{tikzcd}\]
By Corollary~\ref{cor:about-pushout-squares}, all squares of the diagram above are pushout squares. The mapping from multipointed $d$-spaces to multipointed spaces forgetting the set of execution paths being colimit-preserving, being topological, this implies that the map $\overline{X}\to \widehat{Y}$ is a globular subdivision. The composite map $\overline{X}\to \widehat{Y} \to \overline{Y}$ being a globular subdivision by hypothesis, we deduce that the map $S:\widehat{Y}\to \overline{Y}$ is a globular subdivision as well. Thanks to the universal property of the pushout, we obtain the commutative diagram of cellular multipointed $d$-spaces 
\[\begin{tikzcd}[row sep=3em, column sep=3em]
	\globM(\mathbf{S}^{n-1}) \arrow[d,rightarrowtail] \arrow[r,"g"] & X \arrow[d,rightarrowtail] \arrow[r,"\sbd"] & Y \arrow[d,rightarrowtail]\\
	\globM(\mathbf{D}^{n}) \arrow[d,rightarrowtail]\arrow[r,"\widehat{g}"] & \cocartesian \overline{X}\arrow[r,"\sbd"] & \cocartesian\widehat{Y} \arrow[d,rightarrowtail] \arrow[dd,"S",bend left=60pt]\\
	W^n_F \arrow[rr] \arrow[d,"\simeq"']&  \arrow[rd, phantom, "\ulcorner"{font=\Large}, pos=0.9] & T \arrow[d]\cocartesian\\
	\globM(\mathbf{D}^{n})_{F} \arrow[rr] & & \overline{Y} 
\end{tikzcd}\]
By Corollary~\ref{cor:about-pushout-squares}, all squares of the diagram above are pushout squares as well. Using Proposition~\ref{prop:raf-cell}, the map $\globM(\mathbf{S}^{n-1})\to \globM(\mathbf{D}^{n})_{F}$ is a finite composition of pushouts of the maps $C:\varnothing\to\{0\}$ and $\globM(\mathbf{S}^{k-1}) \subset \globM(\mathbf{D}^{k})$ for $k\geq 0$. Thanks to Theorem~\ref{thm:colim}, we obtain the pushout diagram of cellular flows
\[
\begin{tikzcd}[row sep=3em,column sep=3em]
	\glob(\mathbf{S}^{n-1}) \arrow[d,rightarrowtail] \arrow[rr]& \arrow[rd, phantom, "\ulcorner"{font=\Large}, pos=0.9] & \dcat(Y) \arrow[d,rightarrowtail]\\
	\glob(\mathbf{D}^{n})_{F} \arrow[rr] &&  \dcat(\overline{Y})
\end{tikzcd}
\]
By applying Theorem~\ref{thm:colim} again several times, we obtain the composition of pushout diagrams of cellular flows 
\[\begin{tikzcd}[row sep=3em, column sep=4em]
	\glob(\mathbf{S}^{n-1}) \arrow[d,rightarrowtail] \arrow[r,"\dcat(g)"] & \dcat(X) \arrow[d,rightarrowtail] \arrow[r] & \dcat(Y) \arrow[d,rightarrowtail]\\
	\glob(\mathbf{D}^{n}) \arrow[d,rightarrowtail]\arrow[r,"\dcat(\widehat{g})"] & \cocartesian \dcat(\overline{X})\arrow[r] & \cocartesian\dcat(\widehat{Y}) \arrow[d,rightarrowtail] \\
	\dcat(W^n_F) \arrow[rr] &   & \dcat(T) \cocartesian
\end{tikzcd}\]
By composing the latter diagram with the commutative square 
\[
\begin{tikzcd}[row sep=3em,column sep=3em]
	\dcat(W^n_F) \arrow[rr] \arrow[d,"\simeq"']&   & \dcat(T) \arrow[d]\\
	\glob(\mathbf{D}^{n})_{F} \arrow[rr] &  &\dcat(\overline{Y})
\end{tikzcd}
\]
we deduce using Corollary~\ref{cor:about-pushout-squares} the commutative diagram of cellular flows
\[\begin{tikzcd}[row sep=3em, column sep=4em]
	\glob(\mathbf{S}^{n-1}) \arrow[d,rightarrowtail] \arrow[r,"\dcat(g)"] & \dcat(X) \arrow[d,rightarrowtail] \arrow[r] & \dcat(Y) \arrow[d,rightarrowtail]\\
	\glob(\mathbf{D}^{n}) \arrow[d,rightarrowtail]\arrow[r,"\dcat(\widehat{g})"] & \cocartesian \dcat(\overline{X})\arrow[r] & \cocartesian\dcat(\widehat{Y}) \arrow[d,rightarrowtail] \arrow[dd,"\dcat(S)",bend left=60pt]\\
	\dcat(W^n_F) \arrow[rr] \arrow[d,"\simeq"']&  \arrow[rd, phantom, "\ulcorner"{font=\Large}, pos=0.9] & \dcat(T) \arrow[d]\cocartesian\\
	\glob(\mathbf{D}^{n})_{F} \arrow[rr] & & \dcat(\overline{Y}) 
\end{tikzcd}\]
Thanks to the universal property of the pushout, with $V = \dcat(W^n_F)\sqcup_{\glob(\mathbf{D}^n)} U$, we obtain the commutative diagram of cellular flows of Figure~\ref{fig:add-and-subdivide-one-globe} where $\overline{Z}$ is obtained by writing the map  \begin{tikzcd}[cramped,column sep=small]U\arrow[r,rightarrowtail]&V\end{tikzcd} as a composite 
\[
\begin{tikzcd}[column sep=5em,row sep=3em]
	U \arrow[r,rightarrowtail,"\in \cell(\T^{cof})"] & \overline{Z} \arrow[r,"\simeq"]& V
\end{tikzcd}
\]
of a map of $\cell(\T^{cof})$ followed by a weak equivalence of flows thanks to Proposition~\ref{prop:factorization-gl-cof}. The maps $U\to \dcat(\widehat{Y})$ and $V\to \dcat(T)$ of the diagram of Figure~\ref{fig:add-and-subdivide-one-globe} are weak equivalences of the q-model structure of flows since all flows of the diagram above are q-cofibrant and since they are pushouts of the weak equivalence $Z\to \dcat(Y)$ along q-cofibrations (we could also invoke the left properness of the q-model structure of flows \cite[Theorem~5.6]{leftproperflow}). The composite map \[\overline{\imath}:\dcat(\overline{X})\longrightarrow U \longrightarrow \overline{Z}\] belongs to $\cell(\T^{cof})$ since the latter class of maps is closed under composition. It remains to prove that the composite map \[\overline{w}:\overline{Z}\stackrel{\simeq}\longrightarrow V\stackrel{\simeq}\longrightarrow \dcat(T)\longrightarrow \dcat(\overline{Y})\] is a weak equivalence of flows. It therefore remains to prove that the map of flows $\dcat(T)\to \dcat(\overline{Y})$ is a weak equivalence of flows. The commutative square of cellular flows 
\[\begin{tikzcd}[row sep=3em, column sep=3em]
	\dcat(W^n_F) \arrow[d] \arrow[r] & \dcat(T) \arrow[d] \\
	\glob(\mathbf{D}^n)_F \arrow[r] & \cocartesian \dcat(\overline{Y})
\end{tikzcd}\]
being a pushout square, the map $\dcat(T)\to \dcat(\overline{Y})$ induces a bijection on the states, the functor $X\mapsto X^0$ from flows to sets being colimit-preserving and since $\dcat(W^n_F)^0=\glob(\mathbf{D}^n)_F^0=\{0,1\}\cup F$. We obtain the bijections
\[
\dcat(T)^0\iso\dcat(\overline{Y})^0  \iso \dcat(\widehat{Y})^0 \sqcup F.
\]
Consider the set maps $(-)^-,(-)^+:\dcat(T)^0\iso\dcat(\overline{Y})^0 \to \dcat(\widehat{Y})^0$ defined by 
\begin{align*}
	\alpha^-&=\begin{cases}
		\alpha & \hbox{ if }\alpha\in\dcat(\widehat{Y})^0\\
		k(0) & \hbox{ if }\alpha\in F
	\end{cases}\\
	\alpha^+&=\begin{cases}
		\alpha & \hbox{ if }\alpha\in\dcat(\widehat{Y})^0\\
		k(1) & \hbox{ if }\alpha\in F
	\end{cases}
\end{align*}
with $k:\glob(\mathbf{D}^n)\to\dcat(\overline{X})\to U \to \dcat(\widehat{Y})$ (note that it is possible that $k(0)=k(1)$). To complete the proof, it suffices to prove that for all $\alpha,\beta\in T^0$, the induced continuous map $\P_{\alpha,\beta}\dcat(T)\to \P_{\alpha,\beta}\dcat(\overline{Y})$ is a weak homotopy equivalence. Note that it is not possible to invoke the left properness of the q-model category of flows since the map $\dcat(W^n_F)\to \dcat(T)$ has no reason to be a q-cofibration. It then suffices to prove that, for all $\alpha,\beta\in T^0$, the induced continuous map \[\II(\dcat(T))(\alpha,\beta)\to \II(\dcat(\overline{Y}))(\alpha,\beta)\] is a weak homotopy equivalence to complete the proof. Assume that $\alpha,\beta\in F$ and that $\II(\dcat(W^n_F))(\alpha,\beta) \neq \varnothing$. Then one has 
\begin{align*}
	\II(\dcat(T))(\alpha,\beta) & \simeq \{\tau_{\alpha,\beta}\} \sqcup \II(\dcat(\widehat{Y}))(\alpha^+,\beta^-)\\
	& \iso \{\tau_{\alpha,\beta}\} \sqcup \vec{T}(\widehat{Y})(\alpha^+,\beta^-)\\
	& \iso \{\tau_{\alpha,\beta}\} \sqcup \vec{T}(\overline{Y})(\alpha^+,\beta^-) \\
	& \iso \vec{T}(\overline{Y})(\alpha,\beta) \\
	& \iso \II(\dcat(\overline{Y}))(\alpha,\beta),
\end{align*} 
the homotopy equivalence by Theorem~\ref{thm:plus-minus-map-flow}, the first and fourth homeomorphisms by definition of the space of traces, the second homeomorphism by Theorem~\ref{thm:plus-minus-map}, the map $S:\widehat{Y}\longrightarrow\overline{Y}$ being a globular subdivision, and finally the third homeomorphism by Proposition~\ref{prop:recall}~$(3)$. In all other cases, using Theorem~\ref{thm:plus-minus-map-flow} and Theorem~\ref{thm:plus-minus-map} in a similar way, we have 
\begin{align*}
	\II(\dcat(T))(\alpha,\beta) & \simeq \II(\dcat(\widehat{Y}))(\alpha^+,\beta^-)\\
	& \iso  \vec{T}(\widehat{Y})(\alpha^+,\beta^-)\\
	& \iso  \vec{T}(\overline{Y})(\alpha^+,\beta^-) \\
	& \iso  \vec{T}(\overline{Y})(\alpha,\beta) \\
	& \iso \II(\dcat(\overline{Y}))(\alpha,\beta),
\end{align*} 
the third homeomorphism by Proposition~\ref{prop:recall}~$(1,2,4)$.
\epf

\bth \label{thm:facto}
Let \begin{tikzcd}[cramped]f:X\ar[r,"\sbd"]& Y\end{tikzcd} be a globular subdivision of multipointed $d$-spaces. Then $\dcat(f)$ factors as a composite \[\dcat(f)=w.i\] where $i$ belongs to $\cell(\T^{cof})$ and where $w$ is a weak equivalence of the q-model structure of flows. 
\eth

Using the language of \cite{hocont}, the map of flows  $\dcat(f)$ is a dihomotopy equivalence.

\bpf
We consider a cellular decomposition of $X$. Using Theorem~\ref{thm:sbd-cell}, there exists a transfinite tower of cellular multipointed $d$-spaces $\widetilde{Y}:\lambda\to \ptop{\mathcal{M}}$ and a map of transfinite towers $\widetilde{X}\to \widetilde{Y}$ such that the colimit is the globular subdivision $X\to Y$ and such that for all $\nu<\lambda$, there is the 
commutative diagrams of cellular multipointed $d$-spaces of the form
	\[\begin{tikzcd}[row sep=3em, column sep=3em]
	\globM(\mathbf{S}^{n_\nu-1}) \arrow[d,rightarrowtail] \arrow[r,"g_\nu"] \arrow[dd,bend right=60pt,rightarrowtail]& \widetilde{X}_\nu \arrow[d,rightarrowtail] \arrow[r,"\sbd"] & \widetilde{Y}_\nu \arrow[dd,rightarrowtail]\\
	\globM(\mathbf{D}^{n_\nu}) \arrow[d,"\sbd"]\arrow[r,"\widehat{g_\nu}"] & \cocartesian \widetilde{X}_{\nu+1}\arrow[rd,shorten >=0.5em,"\sbd"]  \\
	\globM(\mathbf{D}^{n_\nu})_{F(c_\nu)} \arrow[rr] & & \widetilde{Y}_{\nu+1} \arrow[lu, phantom, "\ulcorner"{font=\Large}, pos=-0.05]
\end{tikzcd}\]
Moreover, the connection maps of the globular subdivision \begin{tikzcd}[cramped]f:\widetilde{X}_\nu\arrow[r,"\sbd"]&\widetilde{Y}_\nu\end{tikzcd} are the restrictions to $Y_\nu^0$ of the connection maps of the globular subdivision \begin{tikzcd}[cramped]f:\widetilde{X}_{\nu+1}\arrow[r,"\sbd"]&\widehat{Y}_{\nu+1}\end{tikzcd} for all $\nu<\lambda$. Let $X_\nu=\dcat(\widetilde{X}_\nu)$ and $Y_\nu=\dcat(\widetilde{Y}_\nu)$. By Theorem~\ref{thm:colim} and Theorem~\ref{thm:rappel-all}~$(4)$, the towers $\nu\mapsto X_\nu$ and $\nu\mapsto Y_\nu$ are two transfinite towers of cellular maps between cellular flows.

We then have to prove by a transfinite induction that, for all $\nu\leq \lambda$, the globular subdivision $\begin{tikzcd}[cramped]\widetilde{X}_\nu\ar[r,"\sbd"]& \widetilde{Y}_\nu\end{tikzcd}$ gives rise to a factorization $X_\nu\to Z_\nu\to Y_\nu$ where the left-hand map belongs to $\cell(\T^{cof})$ and where the right-hand map is a weak equivalence of the q-model structure of flows. The case $\nu=0$ is trivial: $X_0=Y_0=X^0$ indeed. The passage from $\nu<\lambda$ to $\nu+1$ is a consequence of Theorem~\ref{thm:transfinite-sequence-tower}. It remains the case where $\nu\leq\lambda$ is a limit ordinal. The fact that $X_\nu\to Z_\nu$ belongs to $\cell(\T^{cof})$ is a consequence of Theorem~\ref{thm:transfinite-sequence-tower} and Proposition~\ref{prop:hocolim-tower}. The fact that $Z_\nu\to Y_\nu$ is a weak equivalence of the q-model structure of flows is a consequence of Proposition~\ref{prop:homotopy-colimit-tower}.
\epf

\begin{cor} \label{final}
	Let $X$ and $Y$ be two cellular multipointed $d$-spaces related by a finite zigzag sequence of globular subdivisions. Then the associated flows $\dcat(X)$ and $\dcat(Y)$ are related by a finite zigzag of maps of $\cell(\T^{cof})$ and of weak equivalences of flows. 
\end{cor}

Using the language of \cite{hocont}, the two flows  $\dcat(X)$ and $\dcat(Y)$ are dihomotopy equivalent.

\bpf
It is a consequence of Theorem~\ref{thm:facto}.
\epf

Both in Theorem~\ref{thm:facto} and Corollary~\ref{final}, the essentially small class of maps $\T$ can actually be replaced by the essentially small subclass of inclusions of posets $\{0<1\} \subset (\{0,1\} \sqcup F,\leq)$ where $(\{0,1\} \sqcup F,\leq)$ is a globular poset as depicted in Figure~\ref{fig:glposet}, i.e. $F$ is a finite set and the poset structure looks as follows: a smallest element $0$, a biggest element $1$ and a finite number of finite strictly increasing chains of the form $0<a_1<\dots <a_p<1$ with $p\geq 0$.

Intuitively, the underlying homotopy type of a flow is the homotopy type of space obtained after removing the temporal information. It is introduced in \cite[Section~6]{4eme} as follows: 1) we take a q-cofibrant replacement as a cellular flow (there is an arbitrary choice to make); 2) we replace all globes $\glob(\mathbf{D}^n)$ of the cellular decomposition by topological globes $\globM(\mathbf{D}^n)$ attached in the same way (some arbitrary choices have to be made again); 3) we obtain a multipointed $d$-space; 4) we consider the underlying topological space; 5)  it is unique only up to homotopy. The notion is reformulated in \cite[Proposition~8.16]{Moore2}. In \cite{Moore2}, the notion of multipointed $d$-space is slightly different from the one of this paper: the execution paths are reparametrized by nondecreasing homeomorphisms of the segment $[0,1]$, instead of by nondecreasing surjective maps of $[0,1]$. This does not affect the proof of \cite[Proposition~8.16]{Moore2}.

\bd \label{def:underlying-type-flow}
Let $X$ be a flow. The \textit{underlying homotopy type} of $X$ is the homotopy type of the topological space $|(\mathbf{L}\dcat)^{-1}(X)|$ (see the meaning of the notation in Theorem~\ref{thm:rappel-all}~$(1)$). 
\ed

The q-cofibrant multipointed $d$-space $(\mathbf{L}\dcat)^{-1}(X)$ is unique up to weak equivalence. Proposition~\ref{prop:almost-kenbrown} implies that the underlying homotopy type is well-defined.

Proposition~\ref{prop:under-type} is not in \cite{Moore2,Moore3}: it is an omission.

\begin{prop} \label{prop:under-type}
	Let $Z$ be a q-cofibrant multipointed $d$-space, i.e. a retract of a cellular multipointed $d$-space. The underlying homotopy type of the flow $\dcat(Z)$ is the homotopy type of the underlying space $|Z|$ of $Z$.
\end{prop}

\bpf
There are the weak equivalences of multipointed $d$-spaces 
\[
(\mathbf{L}\dcat)^{-1}(\dcat(Z)) \simeq (\mathbf{L}\dcat)^{-1}(\mathbf{L}\dcat)(Z) \simeq Z,
\]
the first weak equivalence since $Z$ is q-cofibrant by hypothesis and by Theorem~\ref{thm:rappel-all}~$(2)$, and the second weak equivalence by Theorem~\ref{thm:rappel-all}~$(1)$. We deduce the homotopy equivalence between q-cofibrant spaces $|(\mathbf{L}\dcat)^{-1}(\dcat(Z))|\simeq |Z|$ by Proposition~\ref{prop:almost-kenbrown}. 
\epf

Proposition~\ref{prop:under-type} is false in general when the q-cofibrancy condition is removed. Indeed, consider the fake crossing $X$ of Figure~\ref{fig:fakecrossing}. Then the homotopy type of $|X|$ is the point whereas the underlying homotopy type of $\dcat(X)=\vI \sqcup \vI$ (see Notation~\ref{nota:glob}) is equal to two points.

We can now use the fact that a globular subdivision is a homeomorphism to prove Corollary~\ref{cor:final2}.

\begin{cor} \label{cor:final2}
	Let $\begin{tikzcd}[cramped]f:X\ar[r,"\sbd"]& Y\end{tikzcd}$ be a globular subdivision. Then the underlying homotopy types of $\dcat(X)$ and $\dcat(Y)$ are equal.
\end{cor}

\bpf
The underlying homotopy type of $\dcat(X)$ ($\dcat(Y)$ resp.) is the homotopy type of $|X|$ ($|Y|$ resp.) by Proposition~\ref{prop:under-type}. The proof is complete since $|f|:|X|\to |Y|$ is a homeomorphism by definition of a globular subdivision.
\epf

\appendix

\section{Ken Brown's lemma}
\label{section:underlying-space}

\cite[Proposition~8.1]{mdtop} wrongly claims that the underlying space functor from multipointed $d$-spaces to topological spaces is a left Quillen functor. In fact it is \textit{almost} a left Quillen functor by Proposition~\ref{prop:underlyingspace-almost-left-Quillen-adjoint} because of the q-cofibration of multipointed $d$-spaces $R:\{0,1\}\to \{0\}$. Proposition~\ref{prop:almost-kenbrown} proves that Ken Brown's lemma is still valid for this almost left Quillen functor.

\begin{prop} \label{prop:double-suspension}
	The endofunctor $Z\mapsto |\globM(Z)|$ of $\top$ is a left Quillen adjoint for the q-model structure of $\top$.
\end{prop}

\bpf 
The calculations are made in the proof of \cite[Theorem~8.2]{4eme} for a different purpose. 
Consider the small category $\C$ (with no composable nonidentity maps)
\[
\begin{tikzcd}[column sep=3em,row sep=0em]
	a & \arrow[l] b \arrow[r,] & c & \arrow[l] d \arrow[r] & e
\end{tikzcd}
\]
equipped with the Reedy structure 
\[
\begin{tikzcd}[column sep=3em,row sep=0em]
	0 & \arrow[l] 1 \arrow[r,] & 2 & \arrow[l] 1 \arrow[r] & 0
\end{tikzcd}
\]
Let us equip the functor category $\top^\C$ of diagrams of spaces over $\C$ with its Reedy q-model structure. If $D$ is an object of the diagram category $\top^\C$, then the latching spaces and the matching spaces of $D$ are as follows: $L_a D= L_b D = L_d D = L_e D = \varnothing$, $L_c D= D_b \sqcup D_d$, $M_a D = M_e D = M_c D = \mathbf{1}$, $M_b D = D_a$ and $M_d D = D_e$. The matching category of an object is either empty or connected. Using \cite[Proposition~15.10.2]{ref_model2} and \cite[Theorem~15.10.8]{ref_model2}, we deduce that the colimit functor $\liminj:\top^\C\to \top$ is a left Quillen adjoint. Consider the functor $\D:\top\to \top^\C$ defined by taking a topological space $Z$ to the diagram of spaces 
\[
\begin{tikzcd}[column sep=3em,row sep=0em]
	\{0\} & \arrow[l] \{0\}\p Z \arrow[r,"\subset"] & {[0,1]}\p Z & \arrow[l,"\supset"'] \{1\}\p Z \arrow[r] & \{1\}
\end{tikzcd}
\]
The functor $\D$ is colimit-preserving, colimits being calculated objectwise in $\top^\C$ and since $\top$ is cartesian closed. Since both the categories $\top$ and $\top^\C$ are locally presentable, it is therefore also a left adjoint by the dual of the Special Adjoint Functor Theorem. A morphism of diagrams $D \longrightarrow E$ of $\top^\C$ is a Reedy q-cofibration when 
\begin{itemize}
	\item $D_a = D_a \sqcup_{L_a D} L_a E \longrightarrow E_a$ is a q-cofibration 
	\item $D_b = D_b \sqcup_{L_b D} L_b E \longrightarrow E_b$ is a q-cofibration 
	\item $D_d = D_d \sqcup_{L_d D} L_d E \longrightarrow E_d$ is a q-cofibration 
	\item $D_e = D_e \sqcup_{L_e D} L_e E \longrightarrow E_e$ is a q-cofibration 
	\item $D_c \sqcup_{(D_b \sqcup D_d)} (E_b\sqcup E_d) = D_c \sqcup_{L_c D} L_c E \longrightarrow E_c$ is a q-cofibration. 
\end{itemize}
Take a q-cofibration $Z_1 \longrightarrow Z_2$ of spaces. The continuous map 
\[([0,1]\p Z_1) \sqcup_{\{0\}\p Z_1 \sqcup \{1\}\p Z_1} (\{0\}\p Z_2 \sqcup \{1\}\p Z_2) \longrightarrow ([0,1]\p Z_2) \] 
is the pushout product of the two q-cofibrations $\{0,1\} \longrightarrow [0,1]$ and $Z_1 \longrightarrow
Z_2$. Thus it is a q-cofibration since the q-model category $\top$ is a monoidal model category by \cite[Proposition~4.2.11]{MR99h:55031}. We deduce that the functor $Z\mapsto \D(Z)$ from $\top$ to $\top^\C$ takes q-cofibrations of spaces to Reedy q-cofibrations. Since it takes any weak homotopy equivalence of spaces to an objectwise weak homotopy equivalence, the functor $\D:\top\to \top^\C$ is therefore a left Quillen adjoint. We deduce that the endofunctor $Z\mapsto |\globM(Z)|$ of $\top$ is a left Quillen adjoint since there is the natural homeomorphism $|\globM(Z)|=\liminj \D(Z)$ for all topological spaces $Z$.
\epf 

\begin{prop} \label{prop:underlyingspace-almost-left-Quillen-adjoint} (replacement for \cite[Proposition~8.1]{mdtop})
	The underlying space functor from multipointed $d$-spaces to topological spaces $X\mapsto |X|$ satisfies the following properties:
	\begin{enumerate}
		\item It is a left adjoint.
		\item It takes q-cofibrations of multipointed $d$-spaces which are one-to-one on states to q-cofibrations of spaces. 
		\item It takes trivial q-cofibrations of multipointed $d$-spaces to trivial q-cofibrations of spaces.
		\item For all q-cofibrant multipointed $d$-spaces $X$, the space $|X|$ is q-cofibrant.
		\item It is not a left Quillen functor.
	\end{enumerate}
\end{prop}

\bpf
The underlying space functor has a right adjoint given by the formula $Z\mapsto (Z,Z,\top([0,1],Z))$. Hence $(1)$. $(2)$ and $(3)$ are consequences of Proposition~\ref{prop:double-suspension}. Since the q-cofibration $R:\{0,1\}\to \{0\}$ is not necessary to reach all cellular multipointed $d$-spaces, we deduce $(4)$. Finally, consider a q-cofibration of multipointed $d$-spaces $f:X\to Y$ which is a nontrivial pushout of the q-cofibration $R:\{0,1\}\to \{0\}$. Then the continuous map $|X|\to |Y|$ is not one-to-one. Therefore it is not a q-cofibration of spaces. Hence $(5)$.
\epf

The point is that the q-cofibration of multipointed $d$-spaces $R:\{0,1\}\to \{0\}$ is not a monomorphism. The latter is not even necessary to reach all cellular objects. So the reader might legitimately ask: why keep it ? Without this (annoying) cofibration, the q-model category of flows collapses by \cite[Corollary~4.10]{Nonunital}. Without it, the minimal model structure has exactly three homotopy types indeed.

\begin{prop} \label{prop:almost-kenbrown} (Ken Brown's lemma for the underlying space functor)
	Let $f:X\to Y$ be a weak equivalence between q-cofibrant multipointed $d$-spaces. Then the continuous map $|f|:|X|\to |Y|$ is a homotopy equivalence between q-cofibrant spaces.
\end{prop}

\bpf We mimic the proof of \cite[Lemma~1.1.12]{MR99h:55031}. Let us start by factoring the map $(f,\id_Y):X\sqcup Y\to Y$ as a composite $\begin{tikzcd}[cramped,column sep=small]X\sqcup Y\arrow[r,rightarrowtail]&Z\arrow[r,twoheadrightarrow,"\simeq"]&Y\end{tikzcd}$ of a q-cofibration followed by a trivial q-fibration. We obtain the commutative diagram of multipointed $d$-spaces 
\[
\begin{tikzcd}[column sep=3em,row sep=3em]
	&X \arrow[rightarrowtail,rd]\arrow[rrd,rightarrowtail,"\simeq"]\arrow[rrrd,"f",bend left=20]\\
	\varnothing \arrow[rr,phantom,pos=0.98,"\scalebox{1.5}{\rotatebox{45}{$\ulcorner$}}"]\arrow[r,rightarrowtail,ru]\arrow[r,rightarrowtail,rd]&& X\sqcup Y  \arrow[r,rightarrowtail] & Z \arrow[r,twoheadrightarrow,"\simeq"]& Y\\
	&Y \arrow[rightarrowtail,ru]\arrow[rru,rightarrowtail,"\simeq"']\arrow[rrru,equal,bend right=20]
\end{tikzcd}
\]
By Proposition~\ref{prop:underlyingspace-almost-left-Quillen-adjoint}, we obtain the commutative diagram of spaces
\[
\begin{tikzcd}[column sep=3em,row sep=3em]
	&{|X|} \arrow[rightarrowtail,rd]\arrow[rrd,rightarrowtail,"\simeq"]\arrow[rrrd,"|f|",bend left=20]\\
	\varnothing \arrow[rr,phantom,pos=0.98,"\scalebox{1.5}{\rotatebox{45}{$\ulcorner$}}"]\arrow[r,rightarrowtail,ru]\arrow[r,rightarrowtail,rd]&& {|X|\sqcup |Y|}  \arrow[r] & {|Z|} \arrow[r,"g"]& {|Y|}\\
	&{|Y|} \arrow[rightarrowtail,ru]\arrow[rru,rightarrowtail,"\simeq"']\arrow[rrru,equal,bend right=20]
\end{tikzcd}
\]
By the \ttt, the map $g$ is a weak homotopy equivalence of spaces. Thus $|f|$ is a weak homotopy equivalence of spaces as well. Since both $|X|$ and $|Y|$ are q-cofibrant spaces by Proposition~\ref{prop:underlyingspace-almost-left-Quillen-adjoint}, the map $|f|$ is a homotopy equivalence.
\epf

Proposition~\ref{prop:almost-kenbrown} is not true in general if the q-cofibrancy condition is removed. Consider the fake crossing $X$ of Figure~\ref{fig:fakecrossing}. Its q-cofibrant replacement is $\vI^{top}\sqcup \vI^{top}$. There is a weak equivalence $f:\vI^{top}\sqcup \vI^{top} \to X$ taking the left copy of $\vI^{top}$ to the directed segment from $(0,1)$ to $(1,0)$ and the right copy of $\vI^{top}$ to the directed segment from $(0,0)$ to $(1,1)$. Thus $|f|$ is a map from a space homotopy equivalent to two points to a contractible space. Therefore $|f|$ is not a weak homotopy equivalence.

%\bibliographystyle{../../plainurlwithoutprefixDOI} 
%\bibliography{../../Bibliotheque}

\end{document}